\documentclass [12pt,a4paper,reqno]{amsart}

\usepackage{geometry}
\usepackage{amsfonts}
\usepackage{amsmath}
\usepackage{amssymb}
\usepackage{amsthm}
\usepackage{amsxtra}
\usepackage{bbm}
\usepackage{braket}
\usepackage{comment}
\usepackage{extarrows}
\usepackage{graphics}
\usepackage{latexsym}
\usepackage{mathrsfs}
\usepackage{mathtools}
\usepackage{verbatim}
\usepackage{nicematrix}
\usepackage[normalem]{ulem}
\usepackage{makecell}
\setcellgapes{3pt}

\usepackage{bm}
\allowdisplaybreaks[4]

\usepackage{tikz}
\usetikzlibrary{patterns}
\usetikzlibrary{decorations.pathreplacing, decorations.markings}
\tikzset{->-/.style={decoration={markings,mark=at position #1 with {\arrow{>}}},postaction={decorate}}}

\usetikzlibrary{arrows}
\usetikzlibrary{patterns,patterns.meta}

\usepackage[initials,lite]{amsrefs} % amsrefs must be loaded after hyperref
\AtBeginDocument{%
    \def\MR#1{}
}
\usepackage[pdfborder={0 0 0}]{hyperref} % hyperref must be loaded after natbib
\hypersetup{colorlinks, citecolor=blue, linkcolor=blue, urlcolor=blue}
\hypersetup{pagebackref}
\usepackage{cleveref}
\Crefname{Lemma}{Lemma}{Lemmas}
\Crefname{Theorem}{Theorem}{Theorems}

\makeatletter
\def\cref@thmoptarg[#1]#2#3#4{%
    \ifhmode\unskip\unskip\par\fi%
    \normalfont%
    \trivlist%
    \let\thmheadnl\relax%
    \let\thm@swap\@gobble%
    \thm@notefont{\fontseries\mddefault\upshape}%
    \thm@headpunct{.}% add period after heading
    \thm@headsep 5\p@ plus\p@ minus\p@\relax%
    \thm@space@setup%
    #2% style overrides
    \@topsep \thm@preskip               % used by thm head
    \@topsepadd \thm@postskip           % used by \@endparenv
    \def\@tempa{#3}\ifx\@empty\@tempa%
      \def\@tempa{\@oparg{\@begintheorem{#4}{}}[]}%
    \else%
      \refstepcounter[#1]{#3}%  <<< cleveref modification
      \@namedef{cref@#3@alias}{#1}% added
      \def\@tempa{\@oparg{\@begintheorem{#4}{\csname the#3\endcsname}}[]}%
    \fi%
    \@tempa}%
\makeatother

 \usepackage{cleveref}
\Crefname{Lemma}{Lemma}{Lemmas}
\Crefname{Theorem}{Theorem}{Theorems}

\makeatletter
\def\cref@thmoptarg[#1]#2#3#4{%
    \ifhmode\unskip\unskip\par\fi%
    \normalfont%
    \trivlist%
    \let\thmheadnl\relax%
    \let\thm@swap\@gobble%
    \thm@notefont{\fontseries\mddefault\upshape}%
    \thm@headpunct{.}% add period after heading
    \thm@headsep 5\p@ plus\p@ minus\p@\relax%
    \thm@space@setup%
    #2% style overrides
    \@topsep \thm@preskip               % used by thm head
    \@topsepadd \thm@postskip           % used by \@endparenv
    \def\@tempa{#3}\ifx\@empty\@tempa%
      \def\@tempa{\@oparg{\@begintheorem{#4}{}}[]}%
    \else%
      \refstepcounter[#1]{#3}%  <<< cleveref modification
      \@namedef{cref@#3@alias}{#1}% added
      \def\@tempa{\@oparg{\@begintheorem{#4}{\csname the#3\endcsname}}[]}%
    \fi%
    \@tempa}%
\makeatother

\usepackage{a4wide}
\theoremstyle{plain}
\newtheorem{Theorem}{Theorem}[section]
\newtheorem{Lemma}[Theorem]{Lemma}
\newtheorem{Corollary}[Theorem]{Corollary}
\newtheorem{Proposition}[Theorem]{Proposition}

\theoremstyle{definition}

\newtheorem{Assumptions and Discussion}[Theorem]{Assumptions and Discussion}

\newtheorem{Example}[Theorem]{Example}
\newtheorem{Definition}[Theorem]{Definition}

\newtheorem{Notation}[Theorem]{Notation}
\newtheorem{remark}[Theorem]{Remark}
\newtheorem{Remark}[Theorem]{Remark}

\theoremstyle{remark}

\newtheorem*{acknowledgment*}{Acknowledgments}

\usepackage[shortlabels]{enumitem}
\SetEnumerateShortLabel{a}{\textup{(\alph*)}}
\SetEnumerateShortLabel{A}{\textup{(\Alph*)}}
\SetEnumerateShortLabel{1}{\textup{(\arabic*)}}
\SetEnumerateShortLabel{i}{\textup{(\roman*)}}
\SetEnumerateShortLabel{I}{\textup{(\Roman*)}}

\definecolor{amaranth}{rgb}{0.9, 0.17, 0.31}
\definecolor{brown(web)}{rgb}{0.65, 0.16, 0.16}
\definecolor{darkcyan}{rgb}{0.0, 0.55, 0.55}
\definecolor{dodgerblue}{rgb}{0.12, 0.56, 1.0}

\def\NN{{\mathbb N}}

\def\N{\operatorname{N}} % nerve complex
\def\N{{\mathbb N}}

\newcommand\bdw{{\bm w}}

\newcommand\bfa{\mathbf{a}}

\newcommand\bfb{\mathbf{b}}

\newcommand\bfv{\mathbf{v}}

\newcommand\bflam{\boldsymbol{\lambda}}

\newcommand\calA{\mathcal{A}}

\newcommand\calC{\mathcal{C}}

\newcommand\calG{\mathcal{G}}
\newcommand\calH{\mathcal{H}}

\usepackage{floatrow}
\usepackage{caption}
\DeclareCaptionSubType[alph]{figure}

\begin{document}


%
%%%%%%%%%%%%%%%%%%%%%%%%%%%%%%%%%%%%%%%%%%%%%%%%%%%%%%%%%%%%%%%%%%%%

\title[IDP of Edge Polytopes of Weighted Oriented Graphs]{Integer Decomposition Property of Edge Polytopes of Weighted Oriented Graphs}

\thanks{Polytope; Integer Decomposition Property; Normal Polytope; Weighted Digraph; Ehrhart ring.}


\author[Biermann]{Jennifer Biermann}
\address{Department of Mathematics and Computer Science, Hobart and William Smith Colleges \\
300 Pulteney St.
Geneva, NY 14456, USA}
\email{biermann@hws.edu}

\author[Lin]{Kuei-Nuan Lin}
\address{Department of Mathematics,
Penn State University\\ 4000 University Dr, McKeesport, PA 15132, USA}
\email{link@psu.edu}

\author[O'Keefe]{Augustine O'Keefe}
\address{Department of Mathematics and Statistics, Connecticut College\\
270 Mohegan Avenue Pkwy,
New London, CT 06320, USA}
\email{aokeefe@conncoll.edu}

\maketitle

\begin{abstract}
We classify the lattice polytopes arising from vertex-weighted oriented graphs that have the integer decomposition property (IDP). Allowing vertex weights and edge orientations introduces non-vertex lattice points in the polytope which need to be accounted for when studying the IDP. We can then identify forbidden induced subgraphs in the sense that their presence obstructs the IDP for the supergraph. We establish positive results for several families of graphs under suitable orientation and coprimality conditions, including weighted oriented paths and cycles, and show that any weighted oriented tree with a unique sink has the IDP.
\end{abstract}

%%%%%%%%%%%%%%%%%%%%%%%%%%%%%%%%%%%%%%%%%%%%%%

\section{Introduction}
In this work, we characterize when a lattice polytope arising from a vertex-weighted oriented graph has the \emph{integer decomposition property} (IDP), a fundamental concept in the study of lattice polytopes. A lattice polytope $P$ is said to have the IDP if every lattice point in a dilation $rP,$ with $r$ a positive integer, can be expressed as a sum of $r$ lattice points from $P$. Although this property is simple to state, determining whether a lattice polytope has the IDP is a challenging problem, as illustrated by the longstanding conjecture by Oda, which asserts that every smooth lattice polytope has the IDP and remains open in general, see e.g.\ \cite{Beck}.

The presence of the IDP has important implications across multiple fields, playing a central role in discrete geometry, optimization, and commutative algebra. In integer programming, a feasible polytope satisfying the IDP simplifies counting arguments and strengthens algorithmic methods for solving combinatorial problems \cite[Chapter 22]{Schrijver1998}. In Ehrhart theory, the presence of the IDP implies that the polytope ring is equal to the Ehrhart ring, which in turn implies that the polytope ring is integrally closed. In the context of toric varieties and commutative algebra, the integral closure of the semigroup ring defined by $P$ ensures desirable algebraic properties such as Cohen-Macaulayness. Within the frameworks of \cite{higashitani2013minkowski} and \cite[Chapter~13]{Sturmfels1996}, the IDP is equivalent to the normality of the semigroup ring.

When polytopes arise from monomials generated by graphs, the connection between the combinatorial structure of the graph and the associated polyhedral geometry becomes especially rich, \cite{OhsugiHibi2003}. In this work we extend the setting to vertex-weighted oriented graphs. In these graphs, each vertex is assigned a weight, and the orientation of edges interacts with these weights to define a monomial algebra. These monomials encode detailed combinatorial information about the graph and simultaneously define lattice polytopes whose vertices correspond to their exponent vectors. Our work also builds on the study of Lin and H\`a, who investigated the IDP of $0$-$1$ polytopes \cite{HL}. In their setting, the generating lattice points are precisely the integer vertices of the polytope. In contrast, our focus requires consideration of non-vertex lattice points. We refer readers to the definitions and notations established in \Cref{Sec:Pre}.

We now outline the main results of this paper. First, we establish the presentation and structure of non-vertex lattice points in \Cref{Sec:NVP} (\Cref{claw}). These results are then used to determine forbidden structures that prevent a polytope associated with a vertex-weighted oriented graph from having the IDP, as described in \Cref{Sec:Forbidden}. This characterization allows us to restrict attention to graphs that do not contain any of the induced subgraphs described in \Cref{SmallestNotNormal}, \Cref{cycleSameOrientation}, and \Cref{LonggcdNot1}.

Next, we present cases in which a polytope of a vertex-weighted oriented graph possesses the IDP. For weighted oriented paths or cycles in which all vertices have the same weight, we show that if the vertices of a path or cycle are either sinks or sources, then the polytope has the IDP (\Cref{PathSinkSource} and \Cref{balanced}). For the case in which vertices have different weights, we apply basic number-theoretic arguments, using coprimality of weights to establish the IDP. In particular, we show that a path or cycle whose vertices are either sinks or sources has the IDP under the additional assumption that the vertex weights are coprime (\Cref{Longgcd1} and \Cref{CycelSinkSource}). Finally, we prove that a weighted oriented tree with a unique sink has the IDP (\Cref{sinkTree}).

\section{Preliminaries}\label{Sec:Pre}

Let $\calG=(V(\calG), E(\calG))$ be a weighted oriented graph with weight vector $\bdw\in  \mathbb{N}^{|V(\calG)|}$ where $V(\calG)=\{x_1,\dots,x_n\}$, $E(\calG)=\{e_1,\dots,e_m\}$, and $\bdw=(w_1,w_2,\dots,w_n)$. We abuse notation and define the polynomial rings $K[E(\calG)]=K[e_1,\dots,e_m]$ and $K[V(\calG)]=K[x_1,\dots,x_n]$ where $K$ is any field. We then consider the $K$-algebra homomorphism $\varphi: K[e_1,\dots,e_m]\rightarrow K[x_1,\dots,x_n]$ defined by $\varphi(e_i)=x_{i_1}x_{i_2}^{w_{i_2}}$ 
where $e_i$ is the directed edge from $x_{i_1}$ to $x_{i_2}$.  The kernel of $\varphi$ is called the \emph{toric ideal} of the weighted oriented graph $\calG$.

\begin{Notation}
    We use the ordered pair $(x_{i_1},x_{i_2})$ to denote the oriented edge from $x_{i_1}$ to $x_{i_2}$. 
    We write $\NN_0=\NN \cup \{0\}$ to be the set of natural numbers with zero.
\end{Notation}

\begin{Definition}
We say that $v\in V(\calG)$ is a \emph{source (or sink)}, if all of the edges incident to it are oriented away from (or toward) $v$.  
\end{Definition}

\begin{Definition}
Let $M(\calG)=(v_{i,j})$ be the $n\times m$ \emph{incidence matrix} of $\calG$ defined by 
\[
    v_{i,j}=
    \begin{cases}
        1 & \text{ if $e_j=(x_i,x_\ell)\in E(\calG)$ for some $1\leq \ell\leq n$},\\
        w_i & \text{ if $e_j=(x_\ell,x_i)\in E(\calG)$ for some $1\leq \ell\leq n$},\\
        0 & \text{ otherwise}
    \end{cases}
\]
with $\{\bfv^1,\dots,\bfv^{m}\}$ the set of columns of $M(\calG)$.  
\end{Definition}

\begin{Notation}
For ease of notation, we switch back and forth between row and column vectors depending on context.
For instance, we denote by $\bfv^i = (v_{i,1}, \ldots, v_{i,n})$ the $i$-th column of $M(\calG)$. Likewise, any vector $\bfa = (a_1, \ldots, a_n) \in \NN^n$ may be regarded as a column vector of a matrix.  In addition, for any vector $\bfa$ we denote the sum of the entries $\bfa$ by $|\bfa|$.
\end{Notation}

\begin{Definition}\label{C_n}
    We call a cycle $\calC_n$ on $n$ vertices \emph{naturally oriented} if all the edges point in the same direction. The incidence matrix of a naturally oriented cycle $\calC_n$ on $n$ vertices (up to a reordering of the vertices) is the following $n\times n$ matrix. 
    \begin{equation}\label{Cn}
       C_n =\left[\begin{array}{cccccc}
 1&w_1 & 0 & \cdots & 0 &0\\
0& 1 & w_2 & \cdots & \vdots &0 \\
0& 0 & 1 & \ddots & 0 &\vdots \\
\vdots & \vdots & \vdots & \ddots &w_{n-2}&0 \\
0& 0 & 0  & \cdots & 1& w_{n-1}\\
w_n & 0 & 0 & \cdots &0& 1 
\end{array}\right].
    \end{equation}
\end{Definition}

We will be looking at the case where the weights are all the same. The following lemma is straightforward to verify, and so we present it without proof. 

\begin{Lemma}\label{invC_n}
    Let $\calC_n$ be a naturally oriented cycle whose vertices all have the same weight $w>1$. If $C_n$ is the incidence matrix of $\calC_n$ as described above, then $\det(C_n)=1+(-1)^{n-1}w^n$, so that it is invertible. Furthermore, its inverse is given by 
    \[
        C_n^{-1}=\frac{1}{\det(C_n)}\begin{bmatrix}
             1 & -w & w^2 & \cdots & (-w)^{n-1}\\
            (-w)^{n-1} & 1 & -w & \cdots & (-w)^{n-2}\\
            \vdots & \vdots & \vdots & \ddots & \vdots\\
            -w & w^2 & -w^3 & \cdots & 1
        \end{bmatrix}.
    \]
\end{Lemma}

Note that the monomials defining the $K$-algebra homomorphism $\varphi$ are exactly the generators of the monomial edge ideal studied in \cites{KBLO, GMSV, HLMRV, MPV, pitones2019monomial}. Defining the toric edge ring of $\calG$ as the image of the monomial map defined by the generators of the monomial edge ideal is a natural analogue to defining the toric edge ring of an unoriented finite simple graph (see \cites{OH, V}).

\begin{remark}\label{rem:disjoint}
If $\calG=\calG_{1}\cup\calG_{2}$ such that $V(\calG_{1})\cap V(\calG_{2})=\emptyset$,
then $I_{\calG}=I_{\calG_1}+ I_{\calG_2}$. Thus we consider graphs with only one connected component.
\end{remark}

We can define a polytope associated to the toric ideal of a weighted oriented graph as follows.

\begin{Definition} \label{normal} 
\begin{enumerate}[a]
    \item Let $\calG$ be a weighted oriented graph with incidence matrix $M(\calG)$ as described above.   We define the \emph{edge polytope} associated to $\calG$ to be the convex hull of the columns of $M(\calG)$.  We denote the edge polytope of $\calG$ by $P_\calG$
    \item A convex polytope $P \subseteq \mathbb{R}^n$ with vertices in the integer lattice $\mathbb{Z}^n$ is said to have \emph{Integer Decomposition Property (the IDP)} if for every positive integer $m$, any point $\mathbf{a}\in mP\cap\mathbf{Z}^n$ can be written as $\mathbf{a}=\sum_{i=1}^m\mathbf{a}_i$ where $\mathbf{a}_i\in P\cap\mathbb{Z}^n$ for all $i=1,2,\dots,m$. 
\end{enumerate}
\end{Definition}

\begin{Remark} 
The IDP is closely related to the normality of a polytope.  The authors in \cite[Chapter 4]{BGPolytopeBook} refer to polytopes having the Integer Decomposition Property as being integrally closed and provide a stricter definition for normality. 
In \cite[Chapter 12]{CCAMillerSturmfelsBook}, having the IDP means a polytope is normal. 
\end{Remark}

The Integer Decomposition Property is a combinatorial property of the polytope $P$ that can be realized algebraically by looking at two rings defined by $P$.

\begin{Definition}  
\begin{enumerate}[a]
    \item The \emph{Ehrhart ring} of a convex lattice polytope is the subalgebra $\mathcal{A}[P] = \bigoplus_{r \geq 0} \mathcal{A}[P]_r$ in the ring $K[x_1,\dots, x_n, z]$ where $\mathcal{A}[P]_r$ is the $K$-vector space spanned by the monomials $\{x^{\bf v}z^r \mid {\bf v} \in rP \cap \mathbb{Z}^n\}$.
    \item The \emph{polytopal ring} of a convex lattice polytope is the monomial subalgebra $K[P] = K[x^{\bf v}z \mid {\bf v} \in P \cap \mathbb{Z}^n]$.
\end{enumerate}
\end{Definition}

\begin{Remark}
It is clear that $K[P]\subseteq \mathcal{A}[P]$. Translating \Cref{normal} into the context of these algebras tells us a convex lattice polytope has the IDP if and only if $\mathcal{A}[P] = K[P]$. 
\end{Remark}

The Integer Decomposition Property of lattice polytopes can be determined using the following criterion. The authors in \cite{HL} worked in the case of 0-1 polytopes, but the proof readily translates to a more general case when considering all lattice points of the polytope.

\begin{Proposition}
\label{Normality}\cite[Proposition 3.1]{HL} Let $P\subseteq\mathbb{R}^{n}$
be a lattice polytope and assume that $E=\{\mathbf{v}^{1},...,\mathbf{v}^{s}\}$
is set of all lattice points of $P$. The polytope $P$ has the IDP if and only
if whenever $\mathbf{a}=\sum_{i=1}^{s}c_{i}\mathbf{v}^{i} \in \mathbb{N}_0^n$ for $0\leq c_{i}<1$
and $\sum_{i=1}^{s}c_{i}\in\mathbb{N}$ then $\mathbf{a}=\sum_{i=1}^{s}d_{i}\mathbf{v}^{i}$
where $d_{i}\in\mathbb{N}_0$ and $\sum_{i=1}^{s}c_{i}=\sum_{i=1}^{s}d_{i}$.
\end{Proposition}

\begin{Notation}
In light of the previous proposition, we will primarily be dealing with exponent vectors in our proofs. We abuse notation by letting $\bfa\in\calA[P]_r$ (respectively $K[P]$) when we mean $x^\bfa z^r\in\calA[P]_r$ (respectively $x^\bfa z \in K[P]$).
\end{Notation}

\section{Non-Vertex Lattice Points}\label{Sec:NVP}

In order to use \Cref{Normality} we need to know not just the vertices of the polytope $P$, but also all of the integer lattice points in $P$. In this section we both characterize when the polytope of a weighted oriented graph $\calG$ has non-vertex lattice points, and in special cases give explicit expressions for them in terms of the vertices of $P_\calG$. 

Note that for any subgraph $\calH$ of $\calG$, a non-vertex lattice point of $P_\calH$ is a non-vertex lattice point of $P_\calG$. Therefore, we focus on finding minimal structures admitting non-vertex lattice points. Our first result describes the general structure of a graph admitting such points.

\begin{Proposition}
      \label{claw}Let $\calG$ be a weighted oriented graph on vertex set $\{x_1, \dots, x_n\}$.  Let $P_\calG$ be the polytope associated to $\calG$ and $\boldsymbol{\lambda}=(\lambda_1,\ldots,\lambda_n)$ be a non-vertex lattice point in $P_\calG$ 
      so that $\boldsymbol{\lambda}=\sum_{i=1}^sc_i\mathbf{v}^i$ where for all $i\in[s]$ we have $\bfv^i$ is a vertex of $P_\calG,$ and $0<c_i<1$ such that $\sum_{i=1}^{s}c_i=1.$ Let $\calH$ be the subgraph of $\calG$ whose incidence matrix is $\begin{bmatrix}\bfv^1&\cdots&\bfv^s\end{bmatrix}.$ Then $\calH$ has at most one source vertex. In the case that there is a source vertex, we can assume, up to a relabeling, that the source vertex is $x_1$ and the edge set of $\calH$ is $\{\bfv^j=(x_1,x_{j+1}):1\leq j\leq s\}.$ 
\end{Proposition}

\begin{proof}
Let $\bflam$ be a non-vertex lattice point of $P_\calG$ as described in the statement of the proposition. Suppose that $\lambda_j$ and $\lambda_{\ell}$ are non-zero and that the corresponding vertices $x_j$ and $x_{\ell}$ are both sources.  Then $\lambda_j = c_{j_1}+\dots+c_{j_a}$ and $\lambda_{\ell} = c_{\ell_1}+\dots+c_{\ell_b}$ with $\{j_1, \dots, j_a\}$ and $\{\ell_1\dots, \ell_b\}$ disjoint.  Since the full sum $\sum_{i = 1}^t c_i$ is equal to one, we have $0\leq \lambda_j, \lambda_{\ell} \leq 1$ and since the two sets of indices are disjoint, they cannot both be equal to 1 which contradicts the assumption that $\lambda_j$ and $\lambda_{\ell}$ are both non-zero integers.

To see the remaining statement, suppose that the subgraph $\calH$ has incidence matrix $\begin{bmatrix}\bfv^1&\cdots&\bfv^s\end{bmatrix}$ as described in the statement of the proposition. Further, assume that $\calH$ has a source vertex that we assume to be $x_1$. Assume the edges incident to $x_1$ in $\calH$ correspond to $\bfv^1,\dots,\bfv^t$ were $t\leq s$. 
Then $\lambda_1 = c_1+\dots +c_{t} \leq c_1+\dots +c_s = 1$. Since $\lambda_1$ is a non-zero integer, we must have $t = s$.
\end{proof}

We now turn our attention to explicit descriptions of the non-vertex lattice points of special subgraphs.  These will be needed in later sections.

\begin{Proposition}\label{NonvertexK12}
Let $\calH$ be the weighted oriented path on three vertices with directed edges $(x_2, x_1)$ and $(x_2, x_3)$ and with weights $w_1$ on vertex $x_1$ and $w_3$ on vertex $x_3$.  Let $\mathbf{v}^1 = (w_1, 1, 0)$ and $\mathbf{v}^2 = (0,1,w_3)$.  Assume that $d = \gcd(w_1, w_3) \neq 1$.  The non-vertex lattice points of $P_\calH$ are exactly the vectors of the form
\[
    \frac{i}{d}\mathbf{v}^1+\frac{d-i}{d}\mathbf{v}^2\text{ for }1\leq i\leq d-1.
\]
\end{Proposition}

\begin{proof}
It is clear that for $1\le i\le d-1$
$$\boldsymbol{\lambda}=\frac{i}{d}\mathbf{v}^1+\frac{d-i}{d}\mathbf{v}^2=\Big(\frac{iw_{1}}{d},1,\frac{(d-i)w_{3}}{d}\Big)\in \mathbb{N}_0^{3}$$ 
and $\frac{i}{d}+\frac{d-i}{d}=1$. Moreover $\boldsymbol{\lambda} \neq \mathbf{v}^{1}, \mathbf{v}^2$. Hence $\boldsymbol{\lambda}$ is a non-vertex lattice point.

Suppose $\boldsymbol{\lambda}=c_1\mathbf{v}^1+c_2\mathbf{v}^2 \in \NN_0^{3}$ is a non-vertex lattice point such that $c_1+c_2=1$ with $0\leq c_1, c_2 <1$. 

 Then $c_{1}w_{1}$ and $c_{2}w_{3}$ are in $\NN_0$ and since $w_{1}, w_{3} \in \N$ this means that $c_{1}$ and $c_{2}$ are rational.  Let $c_{1} = \frac{a}{b}$ be in lowest terms.  Then since $\frac{a}{b}w_{1} \in\NN_0$, we have $a=0$ or $b$ divides $w_{1}$.  Similarly, since $c_{2} = 1-c_{1}$ so that $c_{2}w_3= w_{3} - \frac{a}{b}w_{3} \in \NN_0$, if $a \neq 0$ then $b$ divides $w_{3}$.  Therefore either $a = 0$ or $b$ divides $d = \gcd(w_{1},w_{3})$.  Thus the lattice points of $P_\calH$ are exactly the vectors
 \[
 \boldsymbol{\lambda}=\frac{i}{d}\mathbf{v}^{1}+\frac{d-i}{d}\mathbf{v}^{2}
 \]
 where $0\leq i \leq d$.  Finally, the vertices of $P_\calH$  occur when $i = 0, d$ and therefore $P_\calH$ has non-vertex lattice points if and only if $d \neq 1$. 
\end{proof}

\begin{Corollary}\label{NonvertexPath}
Let $\calG$ be any weighted oriented path.  The non-vertex lattice points of $P_\calG$ are exactly those coming from any subgraph of the form described in \Cref{NonvertexK12}.
\end{Corollary}

\begin{proof}
 Since $\calG$ is a path, every subgraph of $\calG$ has at least one source and the only subgraphs of $\calG$ with exactly one source are a single edge or two edges sharing a source as in \Cref{NonvertexK12}.  Therefore by \Cref{claw} and \Cref{NonvertexK12}, the non-vertex lattice points of $\calG$ are those coming from a subgraph of the form described in \Cref{NonvertexK12}.
\end{proof}

\begin{Example}
    Let $\calH$ be a weighted oriented graph with vertex set $V(\calH)=\{x_1,x_2,\ldots,x_7\}$, edge set $E(\calH)=\{(x_2,x_1),(x_2,x_3),(x_4,x_3),(x_4,x_5),(x_6,x_5),(x_6,x_7)\}$. Further, let weights of the vertices $1,3,5,7$ be $w_1=4,w_3=6,w_5=9,w_7=2$ so that $\gcd(w_1,w_3)=2=d_1,\gcd(w_3,w_5)=3=d_2,\gcd(w_5,w_7)=1=d_3$. 

    \begin{figure}[h]
        \centering
        \includegraphics[scale=1.4]{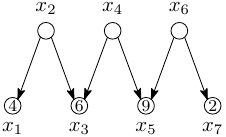}
    \end{figure}
    The vertices of the polytope $P_{\calH}$ are $\mathbf{v}^{1}=(4,1,0,0,0,0,0)$, $\mathbf{v}^{2}=(0,1,6,0,0,0,0)$, $\mathbf{v}^{3}=(0,0,6,1,0,0,0)$, $\mathbf{v}^{4}=(0,0,0,1,9,0,0)$, $\mathbf{v}^{5}=(0,0,0,0,9,1,0)$, $\mathbf{v}^{6}=(0,0,0,0,0,1,2)\in \mathbb{N}_0^{7}$. 
    By \Cref{NonvertexPath}, 
    $$\boldsymbol{\lambda}^1=\Big(2,1,3,0,0,0,0\Big), \boldsymbol{\lambda}^2=\Big(0,0,2,1,6,0,0\Big), \boldsymbol{\lambda}^3=\Big(0,0,4,1,3,0,0\Big)\in \mathbb{N}_0^{7}$$ are non-vertex lattice points of $P_{\calH}$.
\end{Example}

\begin{Lemma}\label{NonVertexClaw}
    Let $\calH$ be the uniformly weighted oriented graph with directed edges $$(x_0, x_1),(x_0,x_2),\ldots,(x_0,x_{n})$$ and weight $w\geq 2$.  Let $\mathbf{v}^i = (1, 0, \ldots, 0, w,0,\ldots,0)$ where the non-zero entries are only at the $0^\text{th}$ and $i^\text{th}$ entries. The non-vertex lattice points of $P_\calH$ are exactly the vectors of the form
\[
    \sum_{i=1}^n \frac{j_i}{w}\mathbf{v}^i=(1,j_1,\ldots, j_n)
\] where $\sum_{i=1}^n {j_i}=w$ and $0\le j_i \le w-1$.
\end{Lemma}

\begin{proof}
    Suppose that $j_i\in\mathbb{N}_0$ with $0\leq j_i\leq w-1$ for $i=1,\dots,n$ such that $\sum_{i=1}^nj_i=w$. It then follows that 
$$\boldsymbol{\lambda}=\sum_{i=1}^n \frac{j_i}{w}\mathbf{v}^i=(1,j_1,\ldots, j_n)\in \mathbb{N}_0^n$$ 
and $\boldsymbol{\lambda} \neq \mathbf{v}^{i}$ for any $i=1,\dots,n$. Hence $\boldsymbol{\lambda}$ is a non-vertex lattice point.

Now suppose $\boldsymbol{\lambda}=\sum_{i=1}^nc_i\mathbf{v}^i \in \NN_0^{n}$ is a non-vertex lattice point such that $\sum_{i=1}^nc_i=1$ with $0\leq c_i <1$. Since $\mathbf{v}^i\in\NN_0^n$ for $i=1,\dots,n$ we must have that the coefficients $c_{i}, $ $i=1,\dots,n$  are rational.  Let $c_{i} = \frac{a_i}{b_i}$ be in lowest terms.  Since $\frac{a_i}{b_i}w \in\NN_0$, we have $a_i=0$ or $b_i$ divides $w$ so that we can write $c_i=\frac{j_i}{w}$ for each $i=1,\dots,n$. Furthermore, since $\sum_{i=1}^n c_i=1$ we must have $\sum_{i=1}^n j_i=w$.  
\end{proof}

\begin{Corollary}\label{NonvertexTree}
Let $\calG$ be any weighted oriented tree whose vertices all have the same weight \(w\geq 2\).  The non-vertex lattice points of $P_\calG$ are exactly those coming from any subgraph of the form described in \Cref{NonVertexClaw}.
\end{Corollary}

\begin{proof}
Let $\calG$ be a weighted oriented tree. Since $\calG$ is a tree, every subgraph of $\calG$ has at least one source and the only subgraphs of $\calG$ with exactly one source are edges sharing a source as in \Cref{NonVertexClaw}.  Therefore by \Cref{claw} and \Cref{NonVertexClaw}, the non-vertex lattice points of $\calG$ are those coming from a subgraph of the form described in \Cref{NonVertexClaw}.
\end{proof}

\begin{Example}
    Let $\calG$ be a weighted oriented graph with vertex set $V(\calG)=\{x_1,x_2,\ldots,x_6\}$ and the edge set $E(\calG)=\{(x_2,x_1),(x_2,x_3),(x_2,x_4),(x_5,x_4),(x_5,x_6)\}$ and all weights equal to 2. The vertices of $P_\calG$ are
    $\mathbf{v}^{1}=(2,1,0,0,0,0)$, $\mathbf{v}^{2}=(0,1,2,0,0,0)$, $\mathbf{v}^{3}=(0,1,0,2,0,0)$, $\mathbf{v}^{4}=(0,0,0,2,1,0)$, $\mathbf{v}^{5}=(0,0,0,0,1,2)\in \mathbb{N}_0^{6}$. 

    \begin{figure}[h]
        \centering
        \includegraphics[scale=1.4]{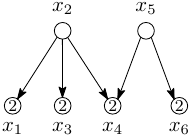}
    \end{figure}
    By \Cref{NonvertexTree}, 
    $$\boldsymbol{\lambda}^1=\Big(1,1,1,0,0,0\Big), \boldsymbol{\lambda}^2=\Big(0,1,1,1,0,0\Big),\boldsymbol{\lambda}^3=\Big(1,1,0,1,0,0\Big),\boldsymbol{\lambda}^4=\Big(0,0,0,1,1,1\Big)\in \mathbb{N}_0^{6}$$ are non-vertex lattice points of $P_{\calG}$. One can see that each induced subgraph has a unique source vertex. 
\end{Example}

For the remainder of this section we consider the case where all vertices in the weighted oriented graph have the same weight. 

\begin{Lemma}\label{w+1}
Let $\calG$ be a weighted oriented graph such that all the vertices of $\calG$ have the same weight $w$. If
$\boldsymbol{\lambda} = (\lambda_1, \lambda_2, \dots, \lambda_n) \in \NN_{0}^n$ is a non-vertex integer lattice point in the polytope $P_{\calG}$ then $|\lambda| = w+1$.
\end{Lemma}
\begin{proof}
Let $\bfv^i$ be the vertex integer lattice points of $P_\calG$ for $1 \leq i \leq t$.
If $\bflam\in \mathbb{N}_0^n$ is a non-vertex lattice point then we can find $0\leq c_i < 1$ such that $\sum_{i=1}^{t}c_i=1$, and $\sum_{i=1}^{t}c_i\bfv^i=\bflam$.  This implies $|\bflam|=\sum_{i=1}^{t}c_i|\bfv^i|=\sum_{i=1}^{t}(c_i(w+1))=(w+1)\sum_{i=1}^{t}c_i=w+1$. 
\end{proof}

\begin{Remark}
    Note that the converse of the previous lemma is not true. Consider the directed graph on three vertices with two edges pointing into a single vertex with weight $2$. It is not possible to write $(1,1,1)$ as a convex combination of the edge vectors. 
\end{Remark}

\begin{Proposition}\label{CyclenonVertexLatticePoints}
    Let $n\geq 3$ and let $\calC_n$ be a cycle such that each vertex has weight $w \geq 2$. Then the non-vertex lattice points of $\calC_n$ either come from the structure described in \Cref{claw} or $\calC_n$ is a naturally oriented cycle and the non-vertex integer lattice points in the polytope $P_{\calC_n}$ are exactly the set
\[
L = \{ \bflam \in \mathbb{N}_0^n \mid \lambda_i \geq 1, |\lambda| = w+1\}.
\]
\end{Proposition}

\begin{proof}
If $\calC_n$ has source vertices, then we can apply \Cref{claw}. 

If $\calC_n$ does not have a source vertex, then $\calC_n$ must be a naturally oriented cycle. By \Cref{w+1}, the non-vertex integer lattice points in the polytope $P_{\calC_n}$ are contained in the set
\[
 \{ \bflam  \in \mathbb{N}_0^n \mid |\lambda| = w+1\}.
\]
Let $\bflam = (\lambda_1, \lambda_2, \dots, \lambda_n)$ be a non-vertex lattice point of $P_{\calC_n}$. If for some $i$, $\lambda_i=0$, the induced subgraph supported on the nonzero entries of $\bflam$ is a disjoint union of paths, each with a single source. By \Cref{claw}, a non-vertex lattice point cannot be supported on such a subgraph, and therefore, $\lambda_i\geq 1$ for all $1\leq i\leq n$. 

It remains to show that every element of the set $L = \{ \bflam \in \mathbb{N}_0^n \mid \lambda_i \geq 1, |\bflam| = w+1\}$ is a non-vertex integer lattice point of $P_{\calC_n}$.  Let $\bflam \in L$.  Let $C_n$ be the incidence matrix of $\calC_n$ as in \Cref{C_n}, and let ${\bf a} = C_n^{-1}\bflam$.  The vector $\bflam$ is a non-vertex integer lattice point of $P_{\calC_n}$ if $|{\bf a}| = 1$ and $a_i \geq 0$ for all $1\leq i\leq n$.

Due to symmetry in the matrix $C_n^{-1}$ described in \Cref{invC_n}, we have that 
\begin{align*}  
    |{\bf a}| = |C_n^{-1}\bflam| 
        &= \frac{1}{1+(-1)^{n-1}w^n} \sum_{i=1}^n \lambda_i (1-w+w^2-\dots+(-w)^{n-1})\\
        &= \frac{1-w+w^2-\dots+(-w)^{n-1}}{1+(-1)^{n-1}w^n} \sum_{i = 1}^n \lambda_i\\
        &= \frac{(1-w+w^2-\dots+(-w)^{n-1})(1+w)}{1+(-1)^{n-1}w^n}\\
        &=1.
\end{align*}

If $n$ is even then $$a_1 = \frac{1}{1-w^n} \sum_{i=1}^{\frac{n}{2}} \left( w^{2i-2}\lambda_{2i-1} - w^{2i-1}\lambda_{2i}\right).$$ Since each $\lambda_i \geq 1$ and $|\bflam| = w+1$, we have that $1\leq \lambda_i \leq w+2-n$ for all $1 \leq i \leq n$.  So $\lambda_{2i-1}w^{2i-2} - \lambda_{2i}w^{2i-1} \leq w^{2i-2}\cdot (w+2-n)-w^{2i-1}\cdot 1= w^{2i-2}(2-n)$.  Since $n\geq 3$ this is negative and $1-w^n$ is also negative, we have $a_1 \geq 0$. 

On the other hand,  if $n$ is odd then $$a_1 = \frac{1}{1+w^n}\left( \lambda_1 + \sum_{i=1}^{\frac{n-1}{2}} \left( -w^{2i-1}\lambda_{2i} + w^{2i}\lambda_{2i+1}\right)\right)$$ and $-w^{2i-1}\lambda_{2i}+w^{2i}\lambda_{2i+1}  \geq -w^{2i-1}(w+2-n) + w^{2i}\cdot 1 =-w^{2i-1}(2-n)$.   Since $n\geq 3$, $-w^{2i-1}(2-n)$ is positive and $1+w^n$ is also positive, we have $a_1 \geq 0$. 

By permuting the elements $\lambda_1 \dots, \lambda_n$, we can see that a similar argument will show that $a_i >0$ for all $i$.  Therefore ${\bflam}$ is in $P_{\calC_n}$.
\end{proof}

\section{Forbidden Structures}\label{Sec:Forbidden}

The next result allows us to find forbidden structures in a weighted oriented graph that prevents its polytope from having the IDP. Namely, if the polytope of an induced subgraph does not have the IDP, neither does the polytope of the whole graph. To leverage this result, we then spend the rest of the section determining minimal structures for which the associated polytope does not have the IDP.

\begin{Theorem}\label{NOTNormal}\cite[Theorem 3.9]{HL}
 Let $\calG$ be a weighted oriented graph and $\calH$ be an induced subgraph of $\calG$. If $P_\calH$ does not have the IDP then $P_\calG$ does not have the IDP. 
\end{Theorem}

\begin{proof}
    Assume $\calH$ is an induced subgraph of $\calG$ and that the polytope $P_\calG$ has the IDP. Since $\calH$ is a subgraph of $P_\calG$, we have that $P_\calH\subseteq P_\calG$. Let $\mathbf{a}$ be a point in the dilation $rP_\calH\subseteq rP_\calG$. Since $P_\calG$ has the IDP, we can write $\mathbf{a}=\sum_{i=1}^r\mathbf{a}_i$ where $\mathbf{a}_i\in P_{\calG}\cap\mathbb{Z}^n$ for all $i=1,2,\dots,r$. In fact, since $\calH$ is an induced subgraph of $\calG$, we can take the sum over $\mathbf{a}_i\in P_{\calH}\cap\mathbb{Z}^n$ for all $i=1,2,\dots,r$ for otherwise, a non-zero entry of $\mathbf{a}$ would correspond to a vertex outside of $\calH$. Thus $P_\calH$ has the IDP.  
\end{proof}

In the next two subsections, we catalogue induced subgraphs that do not have the IDP. 

\subsection{Forbidden Structures with the Same Weight}

In this section we focus on those graphs whose vertices have the same weight $w$.

\begin{Proposition}\label{SmallestNotNormal}
    Let $\calH$ be a weighted oriented graph whose vertices all have the same weight $w>1$ with vertex set 
    $$V(\calH)=\{x_1,x_2,x_3,\ldots,x_{2r+1}\}$$ and edge set 
    \[
    E(\calH)=\{ (x_2,x_1), (x_{2r},x_{2r+1})\}\cup\{(x_{2j+1},x_{2j}), (x_{2j+1},x_{2j+2}) : 1\leq j\leq r-1,\}
    \]
    so that every vertex other than $x_2$ and $x_{2r}$ is either a sink or a source. 
    
    \begin{figure}[h]
\includegraphics[scale=1.4]{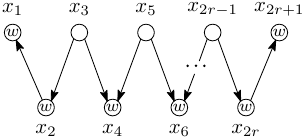}
\end{figure}

Then $\calA[P_\calH]\neq K[P_\calH]$. 
\end{Proposition}

\begin{proof}
 To prove the result, we construct a monomial $x^\bfa z^r \in \calA[P_\calH] \setminus K[P_\calH]$, i.e. $\bfa\in rP\cap \mathbb{Z}^{2r+1}$ but $\bfa\neq \sum_{i=1}^r\bfa^i$ for any collection of integer lattice points $\bfa^i\in P\cap\mathbb{Z}^n$. Then the conclusion will follow from \Cref{Normality}. Let $\bfv^1=(w,1,0,\ldots,0)$, $\bfv^2=(0,w,1,0,\ldots,0)$, $\dots, \bfv^{2r-1}=(0,\ldots,0,1,w,0),$ and $\bfv^{2r}=(0,\ldots,0,1,w)$ be the vertex lattice points of $\calH$. Moreover, by \Cref{NonvertexPath}, the non-vertex lattice points of $P_\calH$ are $$\bflam^{2j+1,k}=(0,\ldots,0,k,1,w-k,0,\ldots,0)$$ for $k=1,\ldots, w-1$ and $j=1,\ldots, r-1$ where the $1$ is in the $(2j+1)$-position. =Let $\bfa=(a_1,\ldots, a_{2r+1})$ where
 \[
    a_j = \begin{cases}
        1, & \text{for }j=2i-1\text{ with }i\in[r],\text{ and }j=2r,\\
        w, & \text{for }j=2i\text{ with }i\in[r-1],\\
        w-1, & \text{for }j=2r+1.
        \end{cases}
\]
 Note that $\bfa$ can be written as a linear combination of the $\bfv^j$'s as follows
    \[
        \bfa =\frac{1}{w}\bfv^1+\frac{w^2-1}{w^2}\sum_{i=1}^{r-1}\bfv^{2i}+\frac{1}{w^2}\sum_{i=1}^{r-1}\bfv^{2i+1}+\frac{w-1}{w}\bfv^{2r}. 
    \]
   Note that $\frac{1}{w}+\sum_{i=1}^{r-1}(\frac{w^2-1}{w^2}+\frac{1}{w^2})+\frac{w-1}{w}=r$ so that $\bfa\in rP\cap\mathbb{Z}^{2r+1}$. 

   Suppose there are $c_i, d_{j,k} \in \mathbb{N}_0$ for $i=1,\ldots,2r$, $j=1,\ldots, r-1$, $k=1,\ldots,w-1$ such that $\sum_{i=1}^{2r}c_i+\sum_{j=1}^{r-1}\sum_{k=1}^{w-1}d_{j,k}=r$ and $\bfa = \sum_{i=1}^{2r}c_{i}\bfv^i+\sum_{j=1}^{r-1}\sum_{k=1}^{w-1}d_{j,k}\bflam^{2j+1,k}$.  Then we will show $c_i=d_{j,k}=0$ for all $i,j,k$ which is a contradiction. Therefore, this implies $x^\bfa z^r \in \calA[P_\calH] \setminus K[P_\calH]$. 
   
   In order to find a solution to the system of equations for $c_i,d_{j,k}$, we consider the augmented matrix 
 \[ A = \begin{bNiceArray}{cw{c}{0.25cm}w{c}{0.25cm}c|c}[margin]
              \Block{2-2}<\Large>{H} & & \Block{2-2}<\Large>{\Lambda}& & \Block{2-1}{\bfa} \\
               & & &  \\ \hline
            1   &\cdots   & \cdots & 1 & r
        \end{bNiceArray}
\]
   where $H$ is the $(2r+1)\times 2r$ incidence matrix of $\calH$ and \[\Lambda=\begin{bmatrix}\bflam^{3,1}&\cdots&\bflam^{3,w-1}&\bflam^{5,1}&\cdots&\bflam^{5,w-1}&\cdots& \bflam^{2r-1,w-1}\end{bmatrix}\] is the $(2r+1)\times(r-1)(w-1)$ matrix whose columns are the nonvertex lattice points. In the submatrix $\begin{bNiceArray}{cc|c} H & \Lambda & \bfa\end{bNiceArray}$, all but the last odd row consist only of ones and zeros. The nonzero entries in these $0$-$1$ rows do not overlap, and adding them together results in the row vector $\begin{bNiceArray}{ccccccccc|c} 0 & 1 & 1  & \cdots & 1& 0 &1 & \cdots &1 &r-1\end{bNiceArray},$ where the second $0$ is in the $2r$-position. Subtracting that from the last row of $A$ results in $c_1+c_{2r}=1$. We must have either $c_1=1$ or $c_{2r}=1$. Suppose $c_1=1$,  then from row $1$ of $A$, we obtain $w=1$, a contradiction. Suppose $c_{2r}=1$, then using row $2r+1$ of $A$, we obtain $c_{2r}w=w-1$, another contradiction. We conclude that $\bfa \in \calA[P_\calH] \setminus K[P_\calH]$.
\end{proof}

\begin{Example}
    To illustrate the proof of \Cref{SmallestNotNormal}, consider the graph on $2r+1=7$ vertices whose structure is as stated in the proposition. Then the matrices $H$ and $\Lambda$ are as follows. 
        \[
  \begin{bNiceArray}{cc}H &\Lambda\end{bNiceArray} = 
            \begin{bNiceArray}{cccccccccccccc}
            \CodeBefore
                \rectanglecolor{blue!15}{3-2}{3-3}
                \rectanglecolor{blue!15}{3-7}{3-10}
                \rectanglecolor{blue!15}{5-4}{5-5}
                \rectanglecolor{blue!15}{5-11}{5-14}
            \Body
             w & 0 & 0 & 0 & 0 & 0 & 0 & 0 & \cdots & 0 & 0 & 0 & \cdots & 0\\
            1 & w & 0 & 0 & 0 & 0 & 1 & 2 & \cdots & w-1 & 0 & 0 & \cdots & 0\\
            0 & 1 & 1 & 0 & 0 & 0  & 1 & 1 & \cdots & 1 & 0 & 0 & \cdots & 0 \\
             0 & 0 & w & w & 0 & 0 & w-1 & w-2 & \cdots & 1 & 1 & 2 & \cdots & w-1 \\
            0 & 0 & 0 & 1 & 1 & 0 & 0 & 0 & \cdots & 0 & 1 & 1 & \cdots & 1 \\
            0 & 0 & 0 & 0 & w & 1 & 0 & 0 & \cdots & 0 & w-1 & w-2 & \cdots & 1\\
            0 & 0 & 0 & 0 & 0 & w & 0 & 0 & \cdots & 0 & 0 & 0 & \cdots & 0
            \end{bNiceArray}
        \]
\end{Example}

\begin{Example}
     Let $\calH$ be a weighted oriented graph with vertex set $V(\calH)=\{x_1,x_2,\ldots,x_5\},$ edge set $E(\calH)=\{(x_2,x_1),(x_3,x_2),(x_3,x_4),(x_4,x_5)\}$, and vertex weights $w_1=3,w_2=3,w_4=3,w_5=3$. Then by \Cref{SmallestNotNormal}, we have $\calA[P_\calH]\neq K[P_\calH]$. If instead, let $\calH'$ be a weighted oriented graph with vertex set $V(\calH')=\{x_1,x_2,\ldots,x_5\},$ edge set $E(\calH')=\{(x_2,x_1),(x_3,x_2),(x_3,x_4),(x_5,x_4)\}$, and vertex weights $w_1=3,w_2=3,w_4=3$. 
     \begin{figure}[h]
\centering
\includegraphics[scale=1.4]{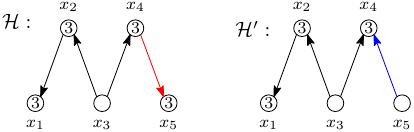}
\end{figure}
     The only difference between $\calH$ and $\calH'$ is the orientation of the edge between vertex $x_4$ and $x_5$. Using Macaulay 2 \cite{M2}, we can verify that  $\calA[P_{\calH'}]= K[P_{\calH'}]$. In fact, one can verify that every other possible orientation of the path graph on five vertices with vertex weight \(3\) has the IDP.
\end{Example}

\begin{Theorem}
\label{cycleSameOrientation} Let $\mathcal{C}_n$ be a naturally oriented $n$-cycle for $n\geq 4$, in which every vertex has weight $w>1$. Then
\[
\mathcal{A}[P_{\calC_n}] \ne K[P_{\calC_n}].
\]
\end{Theorem}

\begin{proof}
By \Cref{CyclenonVertexLatticePoints} the set of non-vertex lattice points of $P_{\calC_n}$ is $\{ \bflam = (\lambda_1, \lambda_2, \dots, \lambda_n) \in \mathbb{N}^n \mid \lambda_i \geq 1, |\lambda| = w+1\}$.  This set is non-empty if and only if $w+1 \geq n$. Suppose $w+1<n$, then there are no non-vertex lattice points. We can
apply the same proof of the case, $w+1 \geq n$.  So suppose $w+1 \geq n$ and there are $s$ elements of this set and label them $\bflam^1, \dots \bflam^s$.

We consider the following vectors:

\[
\mathbf{a}=(\underbrace{1,1,\dots,1}_{\frac{n}{2}},\underbrace{w,\dots,w}_{\frac{n}{2}})
\]
when $n\geq 4$ is even, and 
\[
\mathbf{a}  =(\underbrace{1,1,..,1}_{\lfloor \frac{n}{2}\rfloor},w+1,\underbrace{w...,w}_{\lfloor \frac{n}{2} \rfloor})
\]
when $n\geq 5$ is odd.  
Note that for even $n$, we have $|\mathbf{a}|=\frac{n}{2}(w+1)$ and for odd $n$ we have $|\mathbf{a}|=\left\lceil \frac{n}{2}\right\rceil (w+1)$. We will show that $\mathbf{a}\in \calA[P_{\calC_{n}}]\backslash K[P_{\calC_{n}}]$ then the conclusion follows.

First we show that $\mathbf{a}\in \calA[P_{\calC_{n}}]$, i.e.,  we will show that there is a vector $(c_{1},...,c_{n})$ such that $\mathbf{a}=\sum_{i=1}^{n}c_{i}\mathbf{v}^{i}$
with $0\leq c_{i}<1$ and $\sum_{i=1}^n c_i \in \N$ where $\mathbf{v}^{1}=(1,0,...,0,w), \mathbf{v}^{2}=(w,1,...,0,0),\mathbf{v}^{3}=(0,w,1,0,...,0),\cdots,\mathbf{v}^{n}=(0,0,...,w,1)$.

When $n=2m$, the inverse matrix $C_n$ from \Cref{invC_n} can be rewritten as
\[
        \frac{1}{w^n-1}\begin{bmatrix}
             -1 & w & -w^2 & \cdots & w^{n-1}\\
            w^{n-1} & -1 & w & \cdots & (-1)^{n-1}w^{n-2}\\
            \vdots & \vdots & \vdots & \ddots & \vdots\\
            w & -w^2 & w^3 & \cdots & -1
        \end{bmatrix}.
    \]
    Moreover, one can see that $\mathbf{a}=\sum _{i=1}^n c_i\mathbf{v^i}$ where 
    \begin{equation}\label{evenCycle}
     \begin{bmatrix}
             c_1 \\
            \vdots\\
            c_n
        \end{bmatrix}=\frac{1}{w^n-1}\begin{bmatrix}
             -1 & w & -w^2 & \cdots & w^{n-1}\\
            w^{n-1} & -1 & w & \cdots & (-1)^{n-1}w^{n-2}\\
            \vdots & \vdots & \vdots & \ddots & \vdots\\
            w & -w^2 & w^3 & \cdots & -1
        \end{bmatrix}\begin{bmatrix}
             1 \\
             \vdots\\
            1\\
            w\\
            \vdots\\
            w 
        \end{bmatrix}.
    \end{equation}
    Notice that each $c_i$ is a sum of powers of $w$ with alternating signs. For each $2\leq i\leq n$, the first and last terms cancel, and only for $c_2$ does the largest power of $w$ cancel. In this case the next largest power $(-1)^{n-2}w^{n-2}$ is positive. Therefore each $c_i>0$, since we can pair  terms of opposite sign such that the higher degree term is positive. Finally, one can verify that $\sum_{i=1}^n c_i = m$.

When $n=2m+1$, the inverse matrix $C_n$ from \Cref{invC_n} can be rewritten as
\[
        \frac{1}{w^n+1}\begin{bmatrix}
             1 & -w & w^2 & \cdots & w^{n-1}\\
            w^{n-1} & 1 & -w & \cdots & (-w)^{n-2}\\
            \vdots & \vdots & \vdots & \ddots & \vdots\\
            -w & w^2 & -w^3 & \cdots & 1
        \end{bmatrix}.
    \]
    Moreover, one can see that $\mathbf{a}=\sum _{i=1}^n c_i\mathbf{v^i}$ where 
    \begin{equation}\label{oddCycle}
     \begin{bmatrix}
             c_1 \\
            \vdots\\
            c_n
        \end{bmatrix}=\frac{1}{w^n+1}\begin{bmatrix}
             1 & -w & w^2 & \cdots & w^{n-1}\\
            w^{n-1} & 1 & -w & \cdots & (-w)^{n-2}\\
            \vdots & \vdots & \vdots & \ddots & \vdots\\
            -w & w^2 & -w^3 & \cdots & 1
        \end{bmatrix}\begin{bmatrix}
             1 \\
             \vdots\\
            1\\
            w+1\\
            w\\
            \vdots\\
            w 
        \end{bmatrix}.
    \end{equation}
    Note that each $c_i$ is a sum of powers of $w$ with alternating signs where the largest power of $w$ is positive. We can then see each $c_i$ is positive since for each $i=1,\dots,n$ we have
        \[
            c_i\geq \begin{cases}
                w^{n-1}-\sum_{i=0}^{n-2}w^i & \text{ for }1\leq i\leq \lfloor n/2\rfloor\\
                w^n-\sum_{i=0}^{n-1}w^i & \text{ for }\lceil n/2\rceil\leq i\leq n
            \end{cases}.
        \]
    Furthermore, one can verify that $\sum_{i=1}^n c_i = m+1$.

Now we show that $\mathbf{a}\in \calA[P_{\calC_{n}}]\backslash K[P_{\calC_{n}}]$. Suppose we can rewrite
$\mathbf{a}$ as an integer combinations of the lattice points of $P_{\calC_n}$ so that
\begin{align}
\mathbf{a}
    &=\sum_{i=1}^{n}e_{i}\mathbf{v}^{i}+\sum_{i=1}^{s}d_{i}\boldsymbol{\lambda}^{i} \label{eq:linsys}\\
    & =e_{1}(1,0,...,0,w)+e_{2}(w,1,...,0,0)+\cdots+e_{n}(0,0,...,w,1) \nonumber\\ 
    &\phantom{=} +d_{1}(\lambda_{1}^{1},...,\lambda_{n}^{1})+d_{2}(\lambda_{1}^{2},...,\lambda_{n}^{2})+....+d_{s}(\lambda_{1}^{s},...,\lambda_{n}^{s})\nonumber
\end{align}
where 
\[
    |\mathbf{a}|=\begin{cases}m(w+1), & n=2m\\(m+1)(w+1),& n=2m+1 \end{cases},
\]  
and $e_{i},d_{j}\in\mathbb{N}_0$. Note that \cref{eq:linsys} is a system of linear equations with $n+s$ unknown and $n$
equations. Let $A$ be the system's associated $n\times(n+s+1)$ augmented matrix. 

We show below that $d_j=0$ for all $j$ so that $\mathbf{a}=\sum_{i=1}^{n}e_{i}\mathbf{v}^{i}=\sum_{i=1}^{n}c_{i}\mathbf{v}^{i}$. This gives \(e_2=c_2\) where \(e_2\in\mathbb{N}_0\). But looking at \cref{evenCycle} and \cref{oddCycle} we see that  \(0<c_2<1\), a contradiction.
This shows that $\mathbf{a}\in \calA[P_{\calC_{n}}]\backslash K[P_{\calC_{n}}]$.

We proceed to show that $d_{j}=0$ for all $j$ by proving $e_2<0$ if $d_{j}>0$ for some $j$. Let $C_n$ and $C_n^{-1}$ be as described in \Cref{C_n} and \Cref{invC_n}, and consider the product $C_n^{-1}A$ which yields 
\begin{equation}\label{eq:e_2} 
    e_2 = \frac{1}{\det C_n}\left[-\sum_{j=1}^{s}d_{j}\left(\sum_{i=1}^{n}(-w)^{i-1}\lambda_{i(\text{mod }n)+1}^{j}\right)+B_n\right].
\end{equation}
where 
\[
    B_n = \begin{cases}
        \sum_{i=0}^{m-2}\left[(-w)^{m+i}-(-w)^{i}\right], & n=2m \\ 
        \sum_{i=0}^{m-1}\left[(-w)^{i}-(-w)^{m+i}\right],& n=2m+1 \end{cases}.
\]

Assume that $n=2m+1$ so that $\det C_n>0$. Recall that we are assuming that $d_j\geq 1$ for some $j=1,\dots s$ so that \[-d_j\sum_{i=1}^{n}(-w)^{i-1}\lambda_{i(\text{mod }n)+1}^{j}\leq-\sum_{i=1}^{n}(-w)^{i-1}\lambda_{i(\text{mod }n)+1}^{j}.\]
Thus if for all $j$ we have  
\[
-\sum_{i=1}^{n}(-w)^{i-1}\lambda_{i(\text{mod }n)+1}^{j}+B_n =-\sum_{i=1}^{n}(-w)^{i-1}\lambda_{i(\text{mod }n)+1}^{j} +\sum_{i=0}^{m-1}\left[(-w)^{i}-(-w)^{m+i}\right]<0
\]
we arrive at the desired contradiction, $e_{2}<0$.  In order to show this, we first show that  $-\sum_{i=1}^{n}(-w)^{i-1}\lambda_{i(\text{mod }n)+1}^{j}<0$.
Recall that $1\leq\lambda_i^j\leq w-(n-2)$, and that we are assuming that $w\geq 2$ and $n\geq 4$. We then have
\begin{align*}
-\sum_{i=1}^{n}&(-w)^{i-1}\lambda_{i(\text{mod }n)+1}^{j}\\
&=-\lambda_2^j+w^{n-2}\lambda_{n}^j-w^{n-1}\lambda_1^j+ \sum_{l=2}^{m}\left[w^{2l-3}\lambda_{2l-1}^{j}-w^{2l-2}\lambda_{2l}^{j}\right]\\
&\leq -1 +w^{n-2}(w-(n-2))-w^{n-1}+\sum_{l=2}^{m}\left[w^{2l-3}(w-(n-2))-w^{2l-2}\right]\\
&=-1-\sum_{l=1}^{m}w^{2l-1}(n-2)<0.
\end{align*}
This then tell us 
\begin{align*}
-\sum_{i=1}^{n}(-w)^{i-1}\lambda_{i(\text{mod }n)+1}^{j}+B_n &= -\sum_{i=1}^{n}(-w)^{i-1}\lambda_{i(\text{mod }n)+1}^{j} +\sum_{i=0}^{m-1}\left[(-w)^{i}-(-w)^{m+i}\right]\\
&\leq -1-\sum_{l=1}^{m}w^{2l-1}(n-2)+\sum_{i=0}^{m-1}\left[(-w)^{i}-(-w)^{m+i}\right].\\
\end{align*}

The last expression is a polynomial in $w$ with leading term $-(n-1)w^{2m-1}<0$ yielding the desired contradiction.

When $n$ is even, we have that $\det C_n<0$, and so to arrive at the contradiction $e_2<0$ we aim to show that both 
\[-\sum_{i=1}^{n}(-1)^{i-1}w^{i-1}\lambda_{i(\text{mod }n)+1}^{j}>0,\] and thus
\[-\sum_{i=1}^{n}(-1)^{i-1}w^{i-1}\lambda_{i(\text{mod }n)+1}^{j}+(-1)^{n-1}w^{n-1}+B_n>0.\] The proof is similar to the odd case, and the details are left to the reader.

\end{proof}

\subsection{Forbidden Structure with Different Weights}

\begin{Proposition}\label{LonggcdNot1}

Let $\calH$ be a weighted oriented graph whose edge polytope $P_\calH$ has no non-vertex lattice points with $V(\calH) = \{x_1, \dots, x_{2m+1}\}$ for some $m>1$ and 
\[E(\calH) = \{(x_2, x_1), (x_2, x_3), (x_4, x_3), (x_4, x_5), \dots, (x_{2m}, x_{2m-1}), (x_{2m}, x_{2m+1})\}\] with $\gcd(w_1, w_{2m+1}) \neq 1$.  Then  $\calA[P_\calH]\neq K[P_\calH]$.
\end{Proposition}

\begin{proof}
Let $\calH$ and $V(\calH)$ be as stated above and $E(\calH) = \{e_1= (x_2, x_1), e_2 = (x_2, x_3), e_3 = (x_4, x_3), e_4=(x_4, x_5), \dots, e_{2m-1}=(x_{2m}, x_{2m-1}), e_{2m}=(x_{2m}, x_{2m+1})\}$.  Let ${\bfv^i}$ be the vector associated to the edge $e_i$ for all $1 \leq i \leq t$, and let $d = \mathrm{gcd}(w_1, w_{2m+1})$.  By assumption $d \neq 1$.

We consider the vector,
\begin{align*}
{\bf a} 
    &= \frac{1}{d}{\bfv^1}+ \frac{d-1}{d}{\bfv^2} + \frac{1}{d}{\bfv^3} + \frac{d-1}{d}{\bfv^4} + \dots + \frac{1}{d} {\bfv^{2m-1}}+ \frac{d-1}{d}{\bfv^{2m}}\\ 
    &= \left(\frac{1}{d}w_1, 1, w_3, 1, w_5, 1, \dots, 1, \frac{d-1}{d}w_{2m+1}\right).
\end{align*}
Since $d$ divides both $w_1$ and $w_{2m+1}$, ${\bf a} \in \mathbb{N}^{2m+1}$.  Furthermore, the coefficients of the ${\bfv^i}$ in ${\bfa}$ sum up to $m>1$ giving us $\bfa\in\calA[P_\calH]$.

Now suppose that we can write ${\bfa} = \sum_{i=1}^n c_i{\bfv^i}$ with $c_i \in \mathbb{N}$ and $\sum_{i=1}^n c_i = m$.  There are $m$ entries of ${\bfa}$ equal to one (specifically $a_2, a_4, \dots, a_{2m})$.  Now consider $a_1 = \frac{w_1}{d} \neq 0$. The only edge among $e_1, \dots, e_{2m}$ with $x_1$ as an endpoint is $e_1$, hence $c_1 \neq 0$. But for any integer $n>1$, $n\cdot {\bfv^1}=(nw_1,n,0\ldots,0)$ would give $a_1 = n \cdot w_1 \neq \frac{w_1}{d}$ since $d>1$.   Therefore, it is not possible to write ${\bf a}$ as a positive integer combination of ${\bfv^1}, \dots, {\bfv^t}$ and thus by Proposition \ref{Normality}, $P_\calG$ does not have the IDP.
\end{proof}

\begin{Example}
    Let $\calH$ be a weighted oriented graph with vertex set $V(\calH)=\{x_1,x_2,\ldots,x_7\},$ edge set $E(\calH)=\{(x_2,x_1),(x_2,x_3),(x_4,x_3),(x_4,x_5),(x_6,x_5),(x_6,x_7)\}$, and vertex weights $w_1=4,w_3=5,w_5=3,w_7=2$. We see that $\gcd(w_1,w_3)=\gcd(w_3,w_5)=\gcd(w_5,w_7)=1$, and $\gcd(w_1,w_7)=2$. 

    \begin{center}
        \includegraphics[scale=1.4]{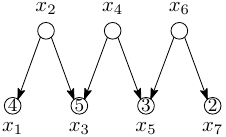}
    \end{center}
    
    By \Cref{NonvertexPath}, 
    there are no non-vertex lattice points of $P_{\calH}$, and the vertices are given by $\mathbf{v}^{1}=(4,1,0,0,0,0,0)$, $\mathbf{v}^{2}=(0,1,5,0,0,0,0)$, $\mathbf{v}^{3}=(0,0,5,1,0,0,0)$, $\mathbf{v}^{4}=(0,0,0,1,3,0,0)$, $\mathbf{v}^{5}=(0,0,0,0,3,1,0)$, $\mathbf{v}^{6}=(0,0,0,0,0,1,2)\in \mathbb{N}^{7}.$ By \Cref{LonggcdNot1},  $\calA[P_\calH]\neq K[P_\calH]$ and its proof shows that \begin{align*}
{\bf a} 
    &= \frac{1}{2}{\bfv^1}+ \frac{1}{2}{\bfv^2} + \frac{1}{2}{\bfv^3} + \frac{1}{2}{\bfv^4} + \frac{1}{2}{\bfv^5} +  \frac{1}{2}{\bfv^{6}}\\ 
    &= \left(2, 1, 5, 1, 3, 1, 1 \right)
\end{align*}
is an integer lattice point of $3P_\calH$ that cannot be written as integer combination of the lattice points of $P_\calH$. 
Note that if instead $w_7=3$ so that $\gcd(w_1,w_7)=1$, then computation with Macaulay2 \cite{M2} shows that the polytope has the IDP. Moreover, it has non-vertex lattice points.
\end{Example}

\section{Polytopes with the IDP}

\subsection{Polytopes of weighted oriented paths}\label{Sec:Normal}

In this section, we present cases of weighted oriented trees that have the IDP. We begin with a simple case.

\begin{Proposition}\label{star}
        Let $\calH$ be a weighted oriented graph such that the vertices of $\calH$ all have the same weight, $w \geq 2$ and such that it has vertex set $V(\calH)=\{x_0,x_1,\ldots,x_{n}\}$, and edge set $E(\calH)=\{(x_0,x_1),(x_0,x_2),\ldots,(x_0,x_{n})\}$. Then $\calA[P_\calH]=K[P_\calH]$.
\end{Proposition}

 \begin{proof}
Let $\mathbf{v}^i$ be the vector associated to the edge $(x_0, x_i)$ for $i\in [n]$. By \Cref{NonVertexClaw} the non-vertex lattice points are exactly the vectors of the form 
 \[
    (1,d_1,d_2,\dots d_n)=\sum_{i=1}^n\frac{d_i}{w}\bfv^i. 
 \]
 where $0\leq d_i \leq w-1$, $ d_i \in \mathbb{N}_{0}$ and $\sum_{i=1}^n d_i=w$. 
 Thus we have  $q={w+n-1 \choose n-1} -n$ non-vertex lattice points. Let $\bflam^j=(1,d_{j1},\ldots,d_{jn})$ for $j\in [q]$ be these non-vertex lattice points. Let $\mathbf{b}=(b_{0},b_{1},\ldots,b_{n})\in \calA[P_{\calH}]$.  We need to show that $\bfb\in K[P_{\calH}]$.
Since $\mathbf{b} \in \calA[P_{\calH}]$ we may write $\mathbf{b}=\sum_{i=1}^{n}c_{i}\mathbf{v}^{i}+\sum_{j=1}^{q}a_j\bflam^{j} \in \NN_{0}^{n+1}$, such that
$\sum_{i=1}^{n}c_{i}+\sum_{j=1}^{q}a_{j} \in \N$, and $0\leq c_{i},a_{j}<1$
for all $i\in[n]$ and $j\in[q]$. By \Cref{Normality}, we need to show $\bfb$ is an integral combination of $\bfv^i$ and $\bflam^{j}$. 

First note that since the first entry in each of $\bfv^i$ and $\bflam^j$ is 1 we have $b_0 = \sum_{i=1}^{n}c_{i}+\sum_{j=1}^{q}a_{j}$.  In addition, since $|\mathbf{v}^i|=w+1$ and  $|\bflam^{j}|=w+1$ for all $i\in [n]$, and $j\in [q]$,  we have $|\bfb| = \sum_{i= 0}^nb_i = b_0(w+1)$ and therefore $\sum_{i = 1}^nb_i = b_0w$.

Define integers $k_{s,t}$ for $1\leq s \leq n$ and $1\leq t \leq b_0$ recursively by $k_{1,1} = \min\{w, b_1\}$ and  $k_{s,t} = \min\{w-\sum_{i = 1}^{s-1}k_{i,t}, b_s - \sum_{j = 1}^{t-1} k_{s,j}\}$ for $1 \leq s \leq n-1$ and $1 \leq t \leq b_0$ and $k_{n,t} = w -\sum_{i = 1}^{n-1}k_{i,t}$ for $1 \leq t \leq b_0$.  Note that  $w-\sum_{i = 1}^{s-1}k_{i,t}\geq 0$ for $2\leq s\leq n$ since otherwise we would have $w-\sum_{i = 1}^{s-2}k_{i,t}<k_{s-1,t}$ which violates the definition of $k_{s-1, t}$.  Similarly, $b_s - \sum_{j = 1}^{t-1} k_{s,j}\geq 0$ for $2\leq t\leq b_0$, and so $k_{s,t} \geq 0$ for all $1\leq s \leq n$ and $1\leq t \leq b_0$.

Note that for any $1 \leq t \leq b_0$, \[\sum_{i = 1}^n k_{i,t} = \left(\sum_{i = 1}^{n-1} k_{i,t}\right) + k_{n,t} = \sum_{i = 1}^{n-1} k_{i,t} + w - \sum_{i =1}^{n-1} k_{i,t} = w.\]  
We claim that $\sum_{j = 1}^{b_0} k_{s, j} = b_s$ for all $1 \leq s \leq n$.  We prove this by induction on $n$.

Let \(s=1\). We first claim that there exists some $1 \leq t_0 \leq b_0$ such that $k_{1,t_0} = b_1 - \sum_{j = 1}^{t_0-1} k_{1,j}$. Suppose that $k_{1,t_0} \neq b_1 - \sum_{j = 1}^{t_0-1} k_{1,j}$ for any $1 \leq t \leq b_0$.  Thus by definition, $k _{1,t} = w$ for all $1 \leq t \leq b_0$, so 
\[
\sum_{j = 1}^{b_0} k_{1,j} = b_0w = \sum_{i =1}^n b_i.
\]
Since $k_{1,b_0} = w = \min\{w, b_1-\sum_{j = 1}^{b_0-1}k_{1,j}\}$ we must have
\[
w < b_1 - \sum_{j = 1}^{b_0-1} k_{1,j}= b_1 - \sum_{j=1}^{b_0-1}w= b_1 - w(b_0 - 1)
\]
which in turn tells us that  \(\sum_{i=1}^n b_i=b_0w < b_1\). This is impossible as each $b_i$ is non-negative integer.

Hence there exists some $1 \leq t_0 \leq b_0$ such that $k_{1,t_0} = b_1 - \sum_{j = 1}^{t_0-1} k_{1,j}$. We have 
\[
k_{1, t_0+1} = \min\left\{w, b_1-\sum_{j = 1}^{t_0}k_{1,j}\right\}
    = \min\left\{w, \left(b_1-\sum_{j = 1}^{t_0-1}k_{1,j}\right)-k_{1,t_0}\right\}
    =0.
\]
A similar argument shows that $k_{1,t} = 0$ for all $t_0 <t\leq b_0$.  Thus
\[
\sum_{j = 1}^{b_0} k_{1,j} = \sum_{j=1}^{t_0}k_{1,j}=b_1.
\] 
This concludes the base case.

Now fix $1<s\leq n$.  Assume that for all $1 \leq i <s$ we have $\sum_{j=1}^{b_0} k_{i,j}= b_i$.  If $s<n$, then we first prove that there exists some $1 \leq t_0 \leq b_0$ such that $k_{s,t_0} = b_s - \sum_{j = 1}^{t_0-1} k_{s,j}$. Suppose that $k_{s,t} \neq b_s - \sum_{j=1}^{t_0-1} k_{s,j}$ for any $1 \leq t \leq b_0$.  Thus $k_{s,t} = w -\sum_{i=1}^{s-1} k_{i,t}$ for all $1 \leq t \leq b_0$.  We then have that
\[
\sum_{j=1}^{b_0} k_{s,j} = \sum_{j=1}^{b_0}\left ( w - \sum_{i=1}^{s-1} k_{i,j}\right )
    = b_0w - \sum_{j=1}^{b_0}\sum_{i=1}^{s-1}k_{i,j}
    = b_0w - \sum_{i=1}^{s-1}\sum_{j=1}^{b_0} k_{i,j}
    = b_0w - \sum_{i = 1}^{s-1} b_i
\]
where the last equality is using the induction hypothesis. Recall that $b_0w=\sum_{i=1}^nb_i,$ which yields
\[
   \sum_{j=1}^{b_0}k_{s,j} = \sum_{i=1}^n b_i - \sum_{i=1}^{s-1} b_i= \sum_{i=s}^n b_i.
\]
Since $k_{s,b_0} = \min\{w-\sum_{i = 1}^{s-1} k_{i, b_0}, b_s - \sum_{j+1}^{b_0-1} k_{s,j}\}$ and $k_{s,b_0} \neq b_s - \sum_{j+1}^{b_0-1} k_{s,j}$ we have
\[
w - \sum_{i=1}^{s-1} k_{i,b_0} < b_s - \sum_{j=1}^{b_0-1} k_{s,j}
    = b_s - \sum_{j=1}^{b_0-1}\left (w - \sum_{i=1}^{s-1}k_{i,j} \right )
    = b_s - w(b_0-1)+\sum_{j=1}^{b_0-1}\sum_{i=1}^{s-1} k_{i,j},
\]

which in turn gives

\[
b_0w < b_s+\sum_{j=1}^{b_0}\sum_{i=1}^{s-1} k_{i,j}
    = b_s + \sum_{i=1}^{s-1}\sum_{j=1}^{b_0} k_{i,j}
    = b_s +\sum_{i=1}^{s-1}b_i = \sum_{i=1}^s b_i.
\]
Therefore $\sum_{i=1}^n b_i = b_0w <\sum_{i=1}^s b_i$ which is impossible as each $b_i$ is a nonnegative integer.

Therefore there exists some $1 \leq t_0 \leq b_0$ such that $k_{s,t_0} = b_s - \sum_{j=1}^{t_0-1} k_{s,j}$. Then 
\[
b_s-\sum_{j = 1}^{t_0}k_{s,j}=\left(b_s-\sum_{j = 1}^{t_0-1}k_{s,j}\right)-k_{s,t_0}=k_{s, t_0}-k_{s,t_0}=0
\]
so that 
\[
k_{s, t_0+1} = \min\left\{w-\sum_{i=1}^{s-1}k_{i, t_0+1}, b_s-\sum_{j = 1}^{t_0}k_{s,j}\right\}=0.
\]
A similar argument shows that $k_{s,t} = 0$ for all $t_0 <t\leq b_0$. Thus we see in this case that 
\[
\sum_{j=1}^{b_0}k_{s,j}=\sum_{j=1}^{t_0}k_{s,j}=b_s.
\]

It remains to show that the summation holds when $s=n$.  By definition $k_{n,t} = w- \sum_{i=1}^{n-1}k_{i,t}$.  Thus, using the induction hypothesis, we have
\[
\sum_{j=1}^{b_0} k_{n,j} = \sum_{j=1}^{b_0} \left (w - \sum_{i=1}^{n-1} k_{i,j}\right )
    = b_0w - \sum_{i=1}^{n-1}\sum_{j=1}^{b_0}k_{i,j}
    = \sum_{i=1}^nb_i - \sum_{i=1}^{n-1}b_i
    = b_n.
\]
Then 
\[
\mathbf{b}  =\left(b_0,\sum_{j=1}^{b_0}k_{1,j},...,\sum_{j=1}^{b_0}k_{n-1,j},\sum_{j=1}^{b_0}\left(w-\sum_{i=1}^{n-1}k_{i,j}\right)\right)
  =\sum_{j=1}^{b_0}\left (1,k_{1,j},\ldots,k_{n-1,j},w-\sum_{i=1}^{n-1}k_{i,j}\right )
\]
is an integral combination of vectors $\mathbf{v}^i$ and $\bflam^j$.  
\end{proof}

\begin{remark}\label{matrices}
Throughout the remainder of the paper, we consider oriented paths whose vertices are all sinks or sources. In the proofs of these results, we are able to reduce to the following structure where $n=2p+1$ for some $p\in\mathbb{N}$. 
\begin{center}
\includegraphics[scale=1.2]{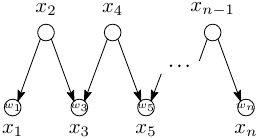}
\end{center}
As such we repeatedly use the associated $n\times(n-1)$ incidence  matrix of the path:
\[
H =\begin{bmatrix}
w_1 & 0 & 0 & \cdots & 0 & 0 \\
1 & 1 & 0 & \cdots  & 0 & 0 \\
0 & w_3 & w_3 & \cdots & 0 & 0 \\
\vdots & \vdots & \vdots & &  \vdots & \vdots \\
0 & 0 & 0 & \cdots & w_{n-2} & 0 \\
0 & 0 & 0 & \cdots &  1 & 1 \\
0 & 0 & 0 & \cdots & 0 & w_n 
\end{bmatrix}.
\]
In the situation where the weights are all $w\in\mathbb{N}$ and there are non-vertex lattice points as described in \Cref{NonvertexK12}, we also consider the following $n\times p(w-1)$ matrix:
\begin{align*}
\Lambda 
&=\begin{bNiceMatrix}[first-row]
{\scriptstyle\bflam^{1,1}} & {\scriptstyle\bflam^{1,2}} & {\scriptstyle\cdots} &{\scriptstyle \bflam^{1,w-1}} & {\scriptstyle\cdots} & {\scriptstyle\bflam^{p,1}} & {\scriptstyle\cdots} & {\scriptstyle\bflam^{p,w-1}}\\
1 & 2 & \cdots & w-1 & \cdots & 0 & \cdots & 0\\
1 & 1 & \cdots & 1 & \cdots & 0 & \cdots & 0 \\
w-1 & w-2 & \cdots & 1 & \cdots & 0 & \cdots & 0 \\
\vdots & \vdots &  & \vdots &  & \vdots &  & \vdots\\
0 & 0 & \cdots & 0 & \cdots & 1 & \cdots & w-1 \\
0 & 0 & \cdots & 0 & \cdots & 1 & \cdots & 1 \\
0 & 0 & \cdots & 0 & \cdots & w-1 & \cdots & 1 
\end{bNiceMatrix}
\end{align*}
\end{remark}

\begin{Proposition}
    \label{PathSinkSource}
    Let $\calG$ be a weighted oriented path whose vertices all have the same weight $w\geq 2$. If all the vertices of $\calG$ are either sinks or sources, then $ \calA[P_{\calG}]=K[P_{\calG}]$.
\end{Proposition}

\begin{proof}
Let $n$ be the number of vertices of $\mathcal{G}$. Denote the vertices of $P_{\mathcal{G}}$ by 
$\mathbf{v}^1, \dots, \mathbf{v}^{n-1}$ , 
and the non-vertex lattice points defined in \Cref{NonvertexPath} by $\boldsymbol{\lambda}^{s,t},$ where $1 \le t \le w-1$ and $1 \le s \le p$, with 
$p \in\left\{ \left\lfloor \frac{n}{2} \right\rfloor,\left\lfloor \frac{n}{2} \right\rfloor - 1\right\}$ being the number of degree two sources in $\mathcal{G}$.

Suppose 
\[\mathbf{a}=(a_1,\ldots,a_n) =\sum_{i=1}^{n-1}c_i\mathbf{v}^i+\sum_{s=1}^{p}\sum_{t=1}^{w-1}d_{s,t}\bflam^{s,t}\in \calA[P_{\calG}]\] 
with $\sum_{i=1}^{n-1}c_i+\sum_{s=1}^{p}\sum_{t=1}^{w-1}d_{s,t}=r\in\mathbb{N}$ such that $r\geq 2$, $0\leq c_i,d_{s,t}<1$ for $1\leq i\leq n-1$, $1\leq s\leq p$, and $1\leq t\leq w-1$. We proceed by rewriting $\mathbf{a}$ as a nonnegative integer combination of the lattice points in $P_{\calG}$ so that by \Cref{Normality} we have $K[P_{\calG}]=\calA[P_{\calG}]$. 

We first make the following reduction. Suppose $\mathbf{v}^1=(1,w,0,\ldots,0)$, then there are no non-vertex lattice points obtained from $\mathbf{v}^1$. Furthermore, since $0\leq c_1<1$, we must have $a_1=c_1=0$. Therefore we may assume $\mathbf{v}^1=(w,1,0,\ldots,0)$. Similarly, we may assume that  $\mathbf{v}^{n-1}=(0,\ldots,0,1,w)$. Furthermore, since every vertex in the path is either a sink or source, we then must have $n=2p+1$ where $p\geq 1$. 

The question is now reduced to showing that there exists a nonnegative integer solution to the system of equations associated to the following augmented matrix: 
\makegapedcells
 \[ A = \begin{bNiceArray}{cw{c}{0.25cm}w{c}{0.25cm}c|c}[margin]
              \Block{2-2}<\Large>{H} & & \Block{2-2}<\Large>{\Lambda}& & \Block{2-1}{\bfa} \\
               & & & \\ \hline
            1   &\cdots   & \cdots & 1 & r
        \end{bNiceArray},
\]
where the matrices $H$ and $\Lambda$ are defined as in \Cref{matrices}.
We perform row operations on $A$ to obtain the following matrix:
\makegapedcells
\[ R= \begin{bNiceArray}{cw{c}{0.25cm}w{c}{0.25cm}c|c}[margin]
              \Block{2-2}<\Large>{H'} & & \Block{2-2}<\Large>{\Lambda'}& & \Block{2-1}{\bfa'} \\
               & & & \\ \hline
            0   & \cdots & \cdots & 0 & r-\sum_{s=1}^{p}a_{2s}
        \end{bNiceArray},
\]
where
        \[
            \begin{bNiceArray}{cc}H' &\Lambda' \end{bNiceArray} = 
            \begin{bNiceArray}{cccccccccccccc}
            \CodeBefore
            \Body
             1 & 0   & \cdots & 0 & 0 & \frac{1}{w} & \frac{2}{w} & \cdots & \frac{w-1}{w} &\cdots & 0 & 0 & \cdots & 0\\
             0 & 1   & \cdots & 0 & 0  & \frac{w-1}{w}  & \frac{w-2}{w}  & \cdots & \frac{1}{w}  & \cdots & 0 & 0 & \cdots & 0 \\
             \vdots & \vdots & \ddots   & \vdots & \vdots & \vdots & \vdots &  & \vdots &  & \vdots & \vdots &  \cdots & 0 \\
             0 & 0   & \cdots & 1 & 0 & 0 & 0 & \cdots & 0 & \cdots &  \frac{1}{w} & \frac{2}{w} & \cdots & \frac{w-1}{w} \\
             0 & 0   & \cdots & 0 & 1 & 0 & 0 & \cdots & 0  & \cdots &  \frac{w-1}{w}  & \frac{w-2}{w} & \cdots & \frac{1}{w} \\
             0 & 0  & \cdots & 0 & 0 & 0 & 0 & \cdots & 0 & \cdots & 0 & 0 & \cdots & 0
            \end{bNiceArray}
        ,\]
 and \(\bfa'=(a'_1,a'_2,\ldots,a'_{2p},a'_{2p+1})\) with
\[
    a'_{2l-1}=\frac{{\textstyle \sum_{j=1}^{l}a_{2j-1}-\sum_{j=1}^{l-1}a_{2j}w}}{w}, ~a'_{2l}=\frac{{\textstyle \sum_{j=1}^{l}a_{2j}w-\sum_{j=1}^{l}a_{2j-1}}}{w}
\]
for all $l\in[p]$ and
\[
a'_{2p+1}=\frac{{\textstyle \sum_{j=1}^{p+1}a_{2j-1}-\sum_{j=1}^{p}a_{2j}w}}{w}.
\]
Since $\bfa\in\mathcal{A}[P_{\mathcal{G}}]$, we know that the system has a solution so that
\[
    a'_{2p+1}=\frac{{\textstyle \sum_{j=1}^{p+1}a_{2j-1}-\sum_{j=1}^{p}a_{2j}w}}{w}=0=r-\sum_{s=1}^{p}a_{2s}. 
\]

From the form of the matrix $R$ we see that the linear system corresponding to rows $(2l-1)$ and $2l$ for each $l\in[p]$ can be written as 
\begin{align}\label{2EqLinSys1}
c_{2l-1} + \frac{1}{w}d_{l,1} + \frac{2}{w}d_{l,2} + \cdots + \frac{w-1}{w}d_{l,w-1} 
    &= a'_{2l-1}\nonumber \\
c_{2l} + \frac{w-1}{w}d_{l,1} + \frac{w-2}{w}d_{l,2} + \cdots + \frac{1}{w}d_{l,w-1} &= a'_{2l}.
\end{align}
Note that the variables in these subsystems do not overlap. Therefore, if we can find nonnegative integer solutions to the system corresponding to rows $(2l-1)$ and $2l$ for each $l\in[p]$, we can combine them to get a nonnegative integer solution to the entire system, thus completing the proof. 

Recall that $0\leq c_i,d_{s,t}<1$ and note that the entries of $\begin{bmatrix}H' & \Lambda'\end{bmatrix}$ are nonnegative. Therefore, by looking at rows $2l-1$ and $2l$ for any $l\in [p]$, we must have 
\[\frac{{\textstyle \sum_{j=1}^{l}a_{2j-1}-\sum_{j=1}^{l-1}a_{2j}w}}{w}, \frac{{\textstyle \sum_{j=1}^{l}a_{2j}w-\sum_{j=1}^{l}a_{2j-1}}}{w} \ge 0.\] 
For each $l\in[p]$ we define 
\[T_l=\sum_{j=1}^{l}a_{2j-1}-\sum_{j=1}^{l-1}a_{2j}w,\] and from the observation immediately above we know that \(T_l \ge 0\). We can then rewrite the system in (\ref{2EqLinSys1}) as 
\begin{align}\label{2EqLinSys2}
c_{2l-1} + \frac{1}{w}d_{l,1} + \frac{2}{w}d_{l,2} + \cdots + \frac{w-1}{w}d_{l,w-1} 
    &= \frac{T_{l}}{w}\nonumber \\
c_{2l} + \frac{w-1}{w}d_{l,1} + \frac{w-2}{w}d_{l,2} + \cdots + \frac{1}{w}d_{l,w-1} &= \frac{a_{2l}w - T_{l}}{w}.
\end{align}
When $T_{l}=kw$ for some $k\in\mathbb{N}_0$ we have the solution 
\[(c_{2l-1}, c_{2l}, d_{l,1}, \dots, d_{l,w-1})=(k,a_{2l}-k,0,\dots,0).\] 
If instead \(T_{l}=k_lw+t_l\) with \(0\le k_l< a_{2l}\) and  $1\le t_l < w$ we have the solution 
\[
(c_{2l-1}, c_{2l}, d_{l,1}, \dots,d_{l,t},\dots, d_{l,w-1})=\left(k_l,a_{2l}-k_l-1,0\dots,1,\dots,0\right),
\]
i.e. $d_{l,t}=1$ and $d_{l,j}=0$ for all $j\neq t$.  
\end{proof}

\begin{Example}
    Let $\calH$ be a weighted oriented graph with vertex set $V(\calH)=\{x_1,x_2,\ldots,x_6\}$ and the edge set $E(\calH)=\{(x_2,x_1),(x_2,x_3),(x_4,x_3),(x_4,x_5),(x_6,x_5)\}$. All vertices have weight 2. Hence $\gcd(w_1,w_3)=\gcd(w_3,w_5)=\gcd(w_1,w_5)=2$. 
    \begin{figure}[h]
    \centering
        \includegraphics[scale=1.4]{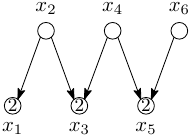}
    \end{figure}
   
   The vertex lattice points of $P_\calH$ are 
    $\mathbf{v}^{1}=(2,1,0,0,0,0)$, $\mathbf{v}^{2}=(0,1,2,0,0,0)$, $\mathbf{v}^{3}=(0,0,2,1,0,0)$, $\mathbf{v}^{4}=(0,0,0,1,2,0)$, $\mathbf{v}^{5}=(0,0,0,0,2,1) \in \mathbb{N}_0^{6}$. By \Cref{claw} and \Cref{NonvertexPath} we have the non-vertex lattice points
    $$\boldsymbol{\lambda}^1=\Big(1,1,1,0,0,0\Big), \boldsymbol{\lambda}^2=\Big(0,0,1,1,1,0\Big)\in \mathbb{N}^{6}.$$   By \Cref{PathSinkSource},  $\calA[P_\calH]= K[P_\calH]$. Notice that \begin{align*}
{\bf a} 
    &= \frac{1}{2}{\bfv^1}+ \frac{1}{2}{\bfv^2} + \frac{1}{2}{\bfv^3} + \frac{1}{2}{\bfv^{4}}\\ 
    &= \left(1, 1, 2, 1, 1, 0\right)=\boldsymbol{\lambda}^1+\boldsymbol{\lambda}^2.
\end{align*}
Note that this graph has an induced subgraph that is almost the structure outlined in \Cref{LonggcdNot1}, e.g., if $w_3=3$ then $\calH$ would not have the IDP. Thus this example shows that the assumption of no non-vertex lattice points is essential in \Cref{LonggcdNot1}. 
\end{Example}

The following number theoretic observation helps when dealing with different weights on the vertices.

\begin{Lemma}\label{coprime}
   Let $w_1, \ldots, w_k \in \mathbb{N}$, and suppose that $w_1, \ldots, w_k$ are pairwise coprime. Suppose \[\sum_{i=1}^k \frac{a_i}{w_i}=r\] for some $a_1,\ldots, a_k,r\in \NN $, then $w_i\mid a_i$ for all $i\in [k]$.
\end{Lemma}
\begin{proof}
   We first multiply $\prod_{i=1}^k w_i$ to the given equation 
and obtain
\[
\sum_{i=1}^k \left (a_i \prod_{j \neq i} w_j\right ) \;=\; r \prod_{i=1}^k w_i.
\]
Notice that for a fixed $i$, the right-hand side is divisible by $w_i$. On the left-hand side, the term $a_p \prod_{j \neq p} w_j$ is divisible by $w_i$ for any $p \neq i$. Therefore, we must have $a_i \prod_{j \neq i} w_j$ is also divisible by $w_i$.  Then by the coprimality assumption $a_i$ is divisible by $w_i$.
\end{proof}

We now present the smallest case in which vertices could have different weights and the polytope satisfies the IDP.

\begin{Lemma} \label{3VNormal}
Let $\calH$ be a weighted oriented graph with $V(\calH)=\{x_1,x_2,x_3\}$ and $E(\calH)=\{(x_2,x_1),(x_2,x_3)\}$. Then $\calA[P_\calH]=K[P_\calH]$.
\end{Lemma}

\begin{proof}
Let $w_1$, $w_3$ be the weights of vertex $1$ and vertex $3$ respectively, and $\mathbf{v}^{1}=(w_1,1,0)$ and $\mathbf{v}^{2}=(0,1,w_3)$
such that $\gcd(w_1,w_3)=d$, then $\bflam^{i}=\frac{i}{d}\mathbf{v}^{1}+\frac{d-i}{d}\mathbf{v}^{2}=(\frac{i}{d}w_1,1,\frac{d-i}{d}w_3)$
are the non-vertex lattice points of $P_\calH$ for all $1\leq i\leq d-1$ by \Cref{NonvertexK12}.  If $d=1$, then there are no non-vertex lattice points. 

Let $\mathbf{b}=(b_{1},b_{2},b_{3})\in \calA[P_{\calH}]$. By \Cref{Normality}, we need to show $\bfb\in K[P_{\calH}]$.
We write $\mathbf{b}=\sum_{i=1}^{2}a_{i}\mathbf{v}^{i}+\sum_{i=1}^{d-1}c_{i}\bflam^i$ with
$\sum_{i=1}^{2}a_{i}+\sum_{i=1}^{d-1}c_{i}=r\in\mathbb{N}$,
and $0\leq a_{i},c_{i}<1$. We consider the augmented matrix
\[
A=\left[\begin{array}{cccccc|c}
w_1 & 0 & \frac{1}{d}w_1 & \frac{2}{d}w_1 & \cdots & \frac{d-1}{d}w_1 & b_{1}\\
1 & 1 & 1 & 1 & \cdots & 1 & b_{2}\\
0 & w_3 & \frac{d-1}{d}w_3 & \frac{d-2}{d}w_3 & \cdots & \frac{1}{d}w_3 & b_{3}\\
1 & 1 & 1 & 1 & \cdots & 1 & r
\end{array}\right].
\]

Since $\mathbf{b}\in\calA[P_\calH]$, subtracting the second row from the fourth row yields $b_2=r$. Further reduction of $A$ yields
\[
A=\left[\begin{array}{ccccc|c}
w_1 & 0 & \frac{1}{d}w_1 & \cdots & \frac{d-1}{d}w_1 & b_{1}\\
1 & 1 & 1 & \cdots & 1 & b_{2}\\
0 & w_3 & \frac{d-1}{d}w_3 & \cdots & \frac{1}{d}w_3 & b_{3}\\
1 & 1 & 1 & \cdots & 1 & r
\end{array}\right]\sim
\left[\begin{array}{ccccc|c}
1 & 1 & 1 & \cdots & 1 & r\\
0 & 1 & \frac{d-1}{d} & \cdots & \frac{1}{d} & r-\frac{b_{1}}{w_1}\\
0 & 0 & 0 & \cdots & 0 & b_{3}+\frac{b_{1}}{w_1}w_3-rw_3\\
0 & 0 & 0 & \cdots & 0 & 0
\end{array}\right]
\]
 \noindent where $\sim$ denotes row equivalence. Therefore $b_{3}+\frac{b_{1}}{w_1}w_3-rw_3=0$, and $\frac{b_3}{w_3}+\frac{b_{1}}{w_1}=r$.
Let $w_1'=w_1/d$ and $w_3'=w_3/d$ then $\gcd(w_1',w_3')=1$ and $\frac{b_3}{w_3'}+\frac{b_{1}}{w_1'}=rd$. Then by \Cref{coprime}, we have  $w_1'|b_{1}$ and $w_3'|b_{3}$, and we can write $b_{1}=k_{1}w_1'$ and $b_{3}=k_{3}w_3'$ for some $k_{1},k_{3}\in\mathbb{N}$. 
This implies that $k_{3}+k_{1}=rd$. Therefore $k_1\le rd$ and we can write $k_{1}=j_{1}+j_{2}+\cdots+j_{r}$
for some $0\leq j_{i}\leq d$ and $j_i\in \mathbb{N}$, which yields $k_{3}=d-j_{1}+d-j_{2}+\cdots+d-j_{r}$. 
Then $\mathbf{b}=(\sum_{i=1}^{r}j_{i}\frac{w_1}{d},r,\sum_{i=1}^{r}(d-j_{i})\frac{w_3}{d})=\sum_{i=1}^{r}(\frac{j_{i}}{d}w_1,1,\frac{d-j_{i}}{d}w_3)=\sum_{i=1}^{r}\bflam^{j_{i}}$
where we take $\bflam^{0}=\bfv^{2}$ and $\bflam^{d}=\mathbf{v}^{1}$. Now we can apply \Cref{Normality} to conclude the proof of the lemma. 
\end{proof}

Notice that when all edges have the same orientation, the three-vertex weighted oriented graph has a unique sink and is a tree; therefore, it has the IDP, as proved later in \Cref{sinkTree}. In the following proposition, we present a special path that has the IDP.

\begin{Proposition}\label{Longgcd1}
Let $\calG$ be the weighted oriented path such that its vertices are either sinks or sources. Moreover, all of its sink vertex weights are pairwise coprime. Then the polytope $P_{\calG}$ has no non-vertex lattice points and $\calA[P_{\calG}]=K[P_{\calG}]$.
\end{Proposition}

\begin{proof}
     By \Cref{NonvertexPath}, the polytope $P_{\calG}$ has no non-vertex lattice points since all of the weights are pairwise coprime.  Let $n$ be the number of vertices of $\calG$ and $\bfv^1,\dots,\bfv^{n-1}$ denote the vertices of $P_{\calG}$.
Suppose $\mathbf{a}=(a_1,\ldots,a_n) =\sum_{i=1}^{n-1}c_i\mathbf{v}^i\in \calA[P_\calG]$ with $\sum_{i=1}^{n-1}c_i=r\in\mathbb{N}$ such that $r\geq 2$, and $0\leq c_i<1$ for $1\leq i\leq n-1$. By \Cref{Normality}, we can show that $K[P_{\calG}]=\calA[P_{\calG}]$ by rewriting $\mathbf{a}$ as a nonnegative integer combination of the lattice points in $P_{\calG}$. 

By the same reduction argument as in the proof of \Cref{PathSinkSource}, we can assume $\mathbf{v}^1=(w_1,1,0,\ldots,0)$, $\mathbf{v}^{n-1}=(0,\ldots,0,1,w_n)$, and $n=2p+1$ for some integer $p\geq 2$. Furthermore, since all vertices are sinks or sources we can assume that $\calG$ has the edge set $E(\calG)=\left\{e_{2l-1}=(x_{2l},x_{2l-1}),e_{2l}=(x_{2l},x_{2l+1})\right\}_{l=1}^{p}.$

Let $\bfa=(a_1,\ldots, a_n)$. Let $H$ be defined as in \Cref{matrices} and consider the following augmented matrix along with its reduced echelon form. 
\makegapedcells
\[           
            =\begin{bNiceArray}{ccc|c}[margin]
              \Block{2-3}<\Large>{H}  && &\Block{2-1}{\bfa} \\
               &  \\ \hline
            1 &  \cdots & 1 & r
        \end{bNiceArray}
        \sim
        \begin{bNiceArray}{ccc|c}[margin]
              \Block{2-3}<\Large>{J} & & &\Block{2-1}{\bfa'} \\
               &  \\ \hline
            0 &  \cdots & 0 & r-\sum_{j = 1}^{\frac{n-1}{2}}a_{2j}
        \end{bNiceArray}
        \]
where $J$ is the block matrix composed of the  $(n-1)\times(n-1)$ identity matrix on top of a row of zeros, and 
\[\bfa'=\left(a_1',a_2',\ldots, a_{2p}',-\sum_{j=1}^{p} a_{2j}+\sum_{j=1}^{p+1} \frac{a_{2j-1}}{w_{2j-1}}\right),\] with 
\[a_{2l-1}'=-\sum_{j=1}^{l-1} a_{2j} +\sum_{j=1}^{l}\frac{a_{2j-1}}{w_{2j-1}}, ~ a_{2l}'=\sum_{j=1}^{l} a_{2j}-\sum_{j=1}^{l}\frac{a_{2j-1}}{w_{2j-1}}\] for $l\in [p]$.

Since the last row of $J$ is comprised of zeros, we must have 
\[
-\sum_{j=1}^{p} a_{2j}+\sum_{j=1}^{p+1} \frac{a_{2j-1}}{w_{2j-1}}=0.
\]
Moreover, since the $c_{i}$'s are nonnegative, we have  $a_{2l-1}',a_{2l}' \ge 0 $ for all $l\in [p]$. Then by \Cref{coprime}, we must have $w_{2j-1}|a_{2j-1}$ when $a_{2j-1}\neq 0$. Therefore, $ a_{2l-1}',a_{2l}' \in \NN_0 $ for all $l\in [p]$. In particular, $\bfa$ can be written as a nonnegative integer combination of the lattice points in $P_{\calG}$.
\end{proof}

Following \Cref{PathSinkSource} and \Cref{Longgcd1}, we present one of the main theorems of this work.

\begin{Theorem}
    Let $\calG$ be a weighted oriented path whose vertices are all either sources or sinks.  Then $P_\calG$ has the IDP if either

\begin{enumerate}[a]
    \item all vertices have the same weight, or
    \item all of the weights on the sinks are relatively prime.
\end{enumerate}
\end{Theorem}

We conclude this subsection with a specific class of trees.

\begin{Theorem}
\label{sinkTree}Let $\calG$ be a weighted oriented tree with a unique
sink and $P_\calG$ denote the associated lattice polytope of $\calG$. Then $\calA[P_\calG]=K[P_\calG]$.
\end{Theorem}

\begin{proof}
By \Cref{claw} and the assumption that $\mathcal{G}$ has a unique sink, we see that $P_{\mathcal{G}}$ has no non-vertex lattice points.

    Let $\{\mathbf{v}^1,\ldots, \mathbf{v}^n\}$ be the set of all lattice points of $P_\calG$. Let $\mathbf{a}=\sum_{i=1}^n c_{i}\mathbf{v}^{i} \in \mathbb{N}_0^{n+1}$ for $0\leq c_{i}<1$
and $\sum_{i=1}^{n}c_{i}\in\mathbb{N}$. Without loss of generality, let $e_1$ be a leaf of $\calG$ with vertex of degree one $x_1$ such that $e_1=(x_1,x_2)$ and  $\mathbf{v}^{1}=(1,w_2,0,\ldots,0)$.  Note that $\mathbf{v}^1$ is the only lattice point among the vertices with a non-zero first coordinate. We must have $a_1=c_1=0$.

Since $\calG$ is a tree with a unique sink, we can pick another leaf of $\calG\setminus \{e_1\}$ and continue this process to conclude that $c_i=0$ for all $i$.
\end{proof}

\subsection{Polytope of weighted oriented cycles}

In the following subsection, we present results related to weighted oriented cycles. We begin with some definitions that we will need.

\begin{Definition}\label{balanced}
\begin{enumerate}[a]
    \item Let $\calC_n$ be a weighted oriented graph. Let $x_i$ be a vertex on $\calC_n$. We say $x_i$ is a clockwise vertex if the adjacent edges of $x_i$ are $(x_{i-1},x_i)$ and $(x_{i},x_{i+1})$. Similarly, we say $x_i$ is a counter-clockwise vertex if the adjacent edges of $x_i$ are $(x_{i},x_{i-1})$ and $(x_{i+1},x_{i})$. 
    \item If $\calC_n$ is a weighted oriented cycle on $n$ vertices and $M(\calC_n)=(v_{i,j})$ is its incidence matrix, we say $\calC_n$ is \emph{balanced} if $n$ is even and $\prod_{i=1}^{n}v_{i,i}=v_{n,1}\prod_{i=1}^{n-1}v_{i,i+1}$.
\end{enumerate}
\end{Definition}

We now present another main theorem of this work.

\begin{Theorem}\label{balanced}
Let $\calC_{n}$ be a weighted oriented balanced cycle such that all vertices have the same weight. Then $P_{\calC_{n}}$ has the IDP if and only if all vertices of $\calC_n$ are either sink or source vertices. 
\end{Theorem}

\begin{proof}
First notice that $n\geq 4$ and $n$ is even. We will make use of Proposition \ref{Normality}.  
We first prove that if all vertices of $\calC_n$ are either sink or source, then $P_{\calC_n}$ has the IDP. 
Without loss of generality, we may assume the vertex vectors of $\calC_{n}$ are 
\[
    \{\mathbf{v}^{1},\mathbf{v}^{2},\dots,\mathbf{v}^{n}\}=\{(w,1,\dots,0,0),(0,1,w,0,\dots,0),\cdots,(0,\dots,0,w,1),(w,\dots,0,0,1)\}.
\]
Let $\frac{n}{2}=p$, and  $\mathbf{a}=(a_1,\ldots,a_n) =\sum_{i=1}^{n}c_{i}\mathbf{v}^{i} +\sum_{s=1}^{p}\sum_{t=1}^{w-1}d_{2s,t}\bflam^{2s,t} \in rP \subseteq \calA[P_{\calC_n}]$ where $\bflam^{2s,t}$ are non-vertex lattices coming from the source vertex $x_{2s}$, and $\sum_{i=1}^{n}c_{i} +\sum_{s=1}^{p}\sum_{t=1}^{w-1}d_{2s,t}=r$ with $0\leq c_i, d_{2s,t}<1$ for all $i\in[n],$ $s\in[p]$, and $t\in[w-1]$. We need to show $\mathbf{a}$ can be written as integer combination of lattice points.

If one of $a_i=0$, then the support of $\mathbf{a}$ is a union of paths such that each path has vertices either sinks or sources. Then by the proof of \Cref{PathSinkSource}, we can show $\mathbf{a}$ is an integer combination of the lattice points of $P_{C_n}$. We now may assume $a_i\geq 1$. We use similar approach in \Cref{PathSinkSource}. 

Let $C=\begin{bmatrix} \bfv^1& \bfv^2 & \cdots & \bfv^n\end{bmatrix}$ be the incidence matrix of $\calC_n$ and $\Lambda$ as described in \Cref{matrices}. We need to show that a system of equations with the following augmented matrix $A$ shown below has non-negative integer solutions: 

\[ A = \begin{bNiceArray}{cw{c}{0.25cm}w{c}{0.25cm}c|c}[margin]
              \Block{2-2}<\Large>{C} & & \Block{2-2}<\Large>{\Lambda}& & \Block{2-1}{\bfa} \\
               & & & \\ \hline
            1   & \cdots  & \cdots & 1 & r
        \end{bNiceArray}.
\]
The matrix $A$ can be reduced to the following matrix $R$:
\makegapedcells
\[ R= \begin{bNiceArray}{cw{c}{0.25cm}w{c}{0.25cm}c|c}[margin]
              \Block{2-2}<\Large>{C'} & & \Block{2-2}<\Large>{\Lambda'}& & \Block{2-1}{\bfa'} \\
               & & & \\ \hline
            0   &\cdots  & \cdots & 0 & r- \sum_{j=1}^{p}a_{2j}
        \end{bNiceArray}
\]
where
\[
    C' = 
    \begin{bNiceArray}{ccccccc}
        \CodeBefore
        \Body
        w & 0  & \cdots & 0 & 0 & 0 & w  \\
        0 & w & \cdots & 0 & 0 & 0  &-w \\
        \vdots & \vdots & \ddots & \vdots & \vdots &\vdots &\vdots\\
        0 & 0 & \cdots & w & 0 &0 &w\\
        0 & 0 &  \cdots & 0 & w &0 &-w \\
        0 & 0  & \cdots & 0 & 0 &w &w\\
        0 & 0  & \cdots & 0 & 0 & 0 &0
    \end{bNiceArray},
\]
\[
    \Lambda'= 
{\scriptsize\begin{bNiceArray}{ccccccccccccccc}
        \CodeBefore
        \Body
         1 & 2 & \cdots & w-1 & 0 & \cdots & 0 & \cdots & 0 & \cdots & 0 & w-1 & w-2 & \cdots & 1\\
         w-1 & w-2 & \cdots & 1 & 0 & \cdots & 0 & \cdots & 0 & \cdots & 0 & 1-w & 2-w & \cdots & -1\\   
         0 & 0 & \cdots & 0 & 1 & \cdots & w-1 & \cdots & 0 & \cdots & 0 & w-1 & w-2 & \cdots & 1\\   
         0 & 0 & \cdots & 0 & w-1 & \cdots & 1 & \cdots & 0 & \cdots & 0 & 1-w & 2-w & \cdots & -1\\   
         \vdots & \vdots &  & \vdots & \vdots &  & \vdots &  & \vdots &  & \vdots & \vdots & \vdots &  & \vdots \\
         0 & 0 & \cdots & 0 & 0 & \cdots & 0 & \cdots & 1 & \cdots & w-1 & w-1 & w-2 & \cdots & 1\\   
         0 & 0 & \cdots & 0 & 0 & \cdots & 0 & \cdots & w-1 & \cdots & 1 & 1-w & 2-w & \cdots & -1\\   
         0 & 0 & \cdots & 0 & 0 & \cdots & 0 & \cdots & 0 & \cdots & 0 & 0 & 0 & \cdots & 0  
    \end{bNiceArray}},
\]
and \(\bfa'=(a'_1,a'_2,\ldots, a'_{2p-1},a'_{2p})\) such that for all $l\in[p]$ we have
\begin{align*}
    a'_{2l-1} &= \sum_{j=1}^la_{2j-1}-w\sum_{j=1}^{l-1}a_{2j},\\
    a'_{2l} &= wa_{2l}-a'_{2l-1}.
\end{align*}
From the $n^\text{th}$ row of $R$ we see that $a'_{2p}=a'_n=0$ so that 
\[
    wa_n = a'_{2p-1} = \sum_{j=1}^p a_{2j-1} - w\sum_{j=1}^{p-1} a_{2j}.
\]
Furthermore, since the linear system associated to these matrices has a nonnegative real solution, namely 
\[
    (c_1,\dots,c_n,d_{2,1},\dots,d_{2,w-1},\dots,d_{p,w-1})\in\mathbb{R}_{\geq0}^{n+p(w-1)},
\] and the entries of each odd row of $\begin{bmatrix}C' & \Lambda'\end{bmatrix}$ are positive, we must also have that $a'_{2l-1}\geq 0$ for all $l\in[p]$. The last nonzero row of $\begin{bmatrix}C' & \Lambda'\end{bmatrix}$ tells us that 
\begin{equation}\label{wan}
w\left(c_{n-1}+c_n+\sum_{j=1}^{w-1}d_{p,j}\right)=wa_n.
\end{equation}

The maximum value of $a_{2l-1}'-wa_{2l}$ will be integral to our construction of a nonnegative solution 
\[
(x_1,\dots,x_n,y_{1,1},\dots,y_{1,w-1}, y_{2,1},\dots,y_{p,w-1})\in\mathbb{N}_0^{n+p(w-1)}
\]
to the linear system. So we define
\[
A = \max\{0, a_1'-wa_2, a_3'-wa_4,\dots,a_{2p-1}'-wa_{2p}\}\geq 0. 
\]
We first show that $A\leq wa_n$. Let $B = \sum_{i=1}^{w-1}(w-i)d_{p,i}$ and note that $B\geq 0$. For each $l\in[p-1]$ the $(2l-1)^\text{th}$ and $(2l)^\text{th}$ rows of $R$ together give us
\[
    a'_{2l-1}-wa_{2l}\leq wc_n +B \leq a'_{2l-1},
\]
and since $wc_n+B\geq 0$, we have for all $l\in[p-1]$, $0\leq A \leq wc_n+B\leq a_{2l-1}'.$ We then have
\begin{equation}\label{x2l}
    wa_{2l}-a'_{2l-1}+A\geq 0.
\end{equation}
Furthermore, since $B\leq (w-1)\sum_{j=1}^{p-1}d_{p,j}$ we see from \Cref{wan} that
\[
    wc_n+B\leq wc_n+w\sum_{j=1}^{p-1}d_{p,j}\leq wa_n
\]
which in turn tells us that
\begin{equation}\label{A_upperBound}
A\leq wa_n.
\end{equation}

We now construct our nonnegative integer solution to the system. From the structure of $R$, we see that we can first determine the values of $x_{n-1}, x_n,y_{p,1},\dots,y_{p,w-1}$, and given those values, we can then treat the system given by the $(2l-1)^\text{th}$ and $(2l)^\text{th}$ rows for each $l\in[p-1]$ independently.

Let $A = wq+s$ with $0\leq s\leq w-1$. If $s=0$, we then define $x_{n-1}=a_n-q$, $x_n=q$, and $y_{p,j}=0$ for all $j\in[w-1]$. Otherwise, define $x_{n-1}=a_n-q-1$, $x_n=q$, and 
\[
    y_{p,j}=
    \begin{cases}
        1 & \text{if }s = w-j\\
        0 & \text{otherwise}
    \end{cases}.
\]
In either case we have 
\[
    x_{n-1}+x_n+\sum_{j=1}^{w-1}y_{p,j}=a_n, 
\]
and
\[
    wx_n+\sum_{j=1}^{w-1}(w-j)y_{p,j}=A.
\]
Furthermore, since $A\leq wa_n$, we must have $q\leq a_n-1$ so we see in particular that $a_n-q-1\geq 0$. Thus, in either case, all of these values are nonnegative. 

For each $l\in[p-1]$, let $a'_{2l-1}-A=wq_l+s_l$ with $0\leq s_l\leq w-1$. Note then that the remainder of $wa_{2l}-a'_{2l-1}+A$ when divided by $w$ is \(0\) if \(s_l=0\), and $w-s_l$ otherwise. If $s_l=0$ we define
\[
    x_{2l-1}=q_l,~~~x_{2l}=a_{2l}-q_l,
\]
and $y_{l,j}=0$ for all $j\in[w-1]$. Otherwise we define
\[
    x_{2l-1}=q_l,~~~x_{2l}=a_{2l}-q_l-1,
\]
and
\[
    y_{l,j}=
    \begin{cases}
        1 & \text{if }j=s_l\\
        0 & \text{otherwise}
    \end{cases}.
\]
Recall from above that $0\leq A\leq a'_{2l-1}$. Along with the inequality in \ref{x2l} we then see that in either case, $x_{2l-1},x_{2l}\geq 0$. Furthermore it is straightforward to check that in any of the four cases regarding the remainders, we have
\begin{align*}
wx_{2l-1}+wx_n+\sum_{j=1}^{w-1}\left[jy_{l,j}+(w-j)y_{p,j}\right] &= a_{2l-1}'\\
wx_{2l}-wx_n+\sum_{j=1}^{w-1}\left[(w-j)y_{l,j}-(w-j)y_{p,j}\right] &= a_{2l}'.
\end{align*}

We now prove the converse of the statement. We prove this by contradiction and assume $\calC_n$ is balanced but there is at least one vertex that is not sink nor source. Since $\calC_n$ is balanced and uniformly weighted, then the number of clockwise and counterclockwise vertices are the same, and all non-sink, non-source vertices cannot be all counter-clockwise or all clockwise, therefore $\calC_n$ has the structure of \Cref{SmallestNotNormal}. This shows that it must be not normal. 
\end{proof}

Similar to \Cref{Longgcd1}, we have the following result for the cycle case.

\begin{Proposition}\label{CycelSinkSource}
Let $\calC_{n}$ be an even cycle such that all vertices are either source
vertices or sink vertices, and all the weights of the sink vertices are pairwise coprime.
Then $K[P_{\calC_{n}}]=\calA[P_{\calC_{n}}]$.
\end{Proposition}

\begin{proof}
Let $n=2p$ with $p>1$. We show $P_{\calC_{n}}$ has no non-vertex lattice points and whenever $\bfa \in \calA[P_{\calC_{n}}]$ then $\bfa \in K[P_{\calC_{n}}]$ at the same time. 

We will make use of Proposition \ref{Normality}. Without loss of generality,
we may assume the vertex vectors of $\calC_{n}$ are \[\{\mathbf{v}^{1},\mathbf{v}^{2},...,\mathbf{v}^{n}\}=\{(w_1,1,...,0),(0,1,w_{3},0,...,0),\ldots,(w_1,...,0,1)\}.\]
Let $\mathbf{a}=\sum_{i=1}^{n}c_{i}\mathbf{v}^{i} \in A[\calC_{n}] $ where $0\leq c_{i}<1$, and $c_{1}+c_{2}+\cdots+c_{n}=r\in\mathbb{N}$. Then 
\begin{align*}
\mathbf{a} & =(a_1,\ldots,a_n)=\sum_{i=1}^{n}c_{i}\mathbf{v}^{i}\\
 & =\left((c_{1}+c_{n})w_1,c_{1}+c_{2},\ldots,(c_{n-2}+c_{n-1})w_{n-1},c_{n-1}+c_{n}\right)\\
\end{align*}
Therefore $c_{2j-1}+c_{2j}=a_{2j}\in\mathbb{N}$, for $j\in [p]$
and $c_{2j}w_{2j+1}+c_{2j+1}w_{2j+1}=a_{2j+1}\in\mathbb{N}$ for $j\in[p-1]$,
and $c_{1}w_{1}+c_{n}w_{1}=a_{1}$. 

Let $C$ be a matrix of the same form as in \Cref{balanced}, but now with different vertex weights, and consider the augmented matrix along with a row equivalent form:
\makegapedcells
\[ A = \begin{bNiceArray}{cc|c}[margin]
              \Block{2-2}<\Large>{C}  & &\Block{2-1}{\bfa} \\
               &  \\ \hline
            1   \cdots & 1 & r
        \end{bNiceArray}\sim \begin{bNiceArray}{cc|c}[margin]
              \Block{2-2}<\Large>{C}  & &\Block{2-1}{\bfa} \\
               &  \\ \hline
            0   \cdots & 0 & r-\sum_{j=1}^{p}a_{2j}
        \end{bNiceArray}=A_1.
\]
In order to have elements in $\calA[P_{\calC_{n}}]$ but
not in $K[P_{\calC_{n}}]$ we must have $r-\sum_{j=1}^{p}a_{2j}=0$. Moreover, $0\leq c_{2j-1}+c_{2j}=a_{2j}<2$ implies $a_{2j}=0$ or $1$. 

We first consider the case that there is some $j\in[p]$ such that $c_{2j}=0$
then $c_{2j-1}+c_{2j}=a_{2j}$ which in turn gives $c_{2j-1}=0=a_{2j}$. By
shifting the indices, we may assume $c_{1}=c_{2}=a_{2}=0$ and $0\leq c_{i}<1$
for $i>2$, then we may consider the matrix $A_{2}$ 
\[ A_2 = \begin{bNiceArray}{cc|c}[margin]
              \Block{2-2}<\Large>{H}  & &\Block{2-1}{\bfa_2} \\
               &  \\ \hline
            0   \cdots & 0 & r-\sum_{j=1}^{p}a_{2j}
        \end{bNiceArray}
\]
where $H$ is obtained by deleting the second row and first two columns of $C$, and $\bfa_2=(a_1, a_3,\ldots,a_n)$.
Note that $H$ is now the incidence matrix for a path, and we now can apply \Cref{Longgcd1} to conclude that $P_{\calC_n}$ has no non-vertex lattice points and $A[P_{\calC_n}]=K[P_{\calC_n}]$.

We now consider the case that $c_{i}\neq0$ for all $i$, then $a_{i}\neq0$
for all $i$. Like before, since $c_{i}\neq0$ for all $i$, we must
have $c_{2j-1}+c_{2j}=a_{2j}=1$ for $j=1,...,p$ and $r=p$.
Then $A_{1}$ is row equivalent to

\[ A_1' = \begin{bNiceArray}{cc|c}[margin]
              \Block{2-2}<\Large>{C'}  & &\Block{2-1}{\bfa'} \\
               &  \\ \hline
            0   \cdots & 0 & 0
        \end{bNiceArray}
        = \begin{bNiceArray}{cccccc|c}[margin]
        1 & 0 & 0 & \cdots & 0 & 1 & a'_1\\
        0 & 1 & 0 & \cdots & 0 & -1 & a'_2\\
        0 & 0 & 1 & \cdots & 0 & 1 & a'_3\\
        \vdots & \vdots & \vdots & \ddots & \vdots & \vdots & \vdots\\
        0 & 0 & 0 & \cdots & 1 & 1 & a'_{2p-1}\\
        0 & 0 & 0 & \cdots & 0 & 0 & a'_{2p}\\ \hline
        0 & 0 & 0 & \cdots & 0 & 0 & 0
        \end{bNiceArray}
\] 
and
        \[
            \bfa'=\left(a'_1,a'_2,\ldots,a'_{2p-3},a'_{2p-2}, a'_{2p-1}, \sum_{j=1}^p\left(a_{2j}-\frac{a_{2j-1}}{w_{2j-1}}\right) \right),
        \]
        such that for all $l\in[p]$
        \[
            a'_{2l-1}= -(l-1)+\sum_{j=1}^{l}  \frac{a_{2j-1}}{w_{2j-1}},~~
            a'_{2l}=l-\sum_{j=1}^{l}  \frac{a_{2j-1}}{w_{2j-1}}.
        \]
Since the last row of $C'$ is comprised of zeros, we see that
\[
\sum_{j=1}^p\left(a_{2j}-\frac{a_{2j-1}}{w_{2j-1}}\right)=p-\sum_{j=1}^p\frac{a_{2j-1}}{w_{2j-1}}=0,
\]
which tells us that $a'_{2p-1}=1$.

By \Cref{coprime} and our assumption that the weights are coprime, this gives $w_{2j-1}|a_{2j-1}$ for all $j\in[p]$. From the row reduced matrix, we have $c_{2p-1}=1-c_{2p},$ and for all $l\in[p-1]$ 
\begin{align*}
c_{2l-1}&=-c_{2p} -(l-1)+\sum_{j=1}^{l}  \frac{a_{2j-1}}{w_{2j-1}},\\
c_{2l}&=c_{2p}+l-\sum_{j=1}^{l}  \frac{a_{2j-1}}{w_{2j-1}}.
\end{align*} 

If $a_{1}/w_{1}\geq2$, then $c_{2}=c_{2p}+1-a_{1}/w_{1}<0$,
a contradiction. Hence we must have $a_{1}/w_{1}=1$, therefore
$c_{2}=c_{2p}$ and $c_{1}=-c_{2p}+1$. 
We iterate this process to see that $a_{2j-1}/w_{2j-1}=1$ for all $j\in [p-1]$ and obtain $c_{2j}=c_{2p}$ and $c_{2j-1}=-c_{2p}+1$.
With $p-\sum_{j=1}^p\frac{a_{2j-1}}{w_{2j-1}}=0$, we conclude that $a_{2j-1}/w_{2j-1}=1$ for all $j\in [p]$. We then have $\bfa=(w_1,1,w_3,\ldots, w_{2p-1},1),$ and we can write
 
\begin{align*}
\mathbf{a} 
 & =(w_1,1,w_3,\ldots,w_{2p-1},1)\\
 & =\sum_{l=1}^{p}\mathbf{v}^{2l-1}.
\end{align*}
This shows that there can't be any non-vertex lattice points, and that any element in $\calA[P_{\calC_{n}}]$
can be written as integer combinations of vertex lattice points of $P_{\calC_{n}}$,
i.e. $K[P_{\calC_{n}}]=\calA[P_{\calC_{n}}]$. 
\end{proof}

We conclude this work with a cycle example.

\begin{Example}
    Let $\calC_6$ be an even cycle weighted oriented graph with vertex set $V(\calC_6)=\{x_1,x_2,\ldots,x_6\}$ and the edge set $E(\calC_6)=\{(x_2,x_1),(x_2,x_3),(x_4,x_3),(x_4,x_5),(x_6,x_5),(x_6,x_1)\}$. The weight of vertices are \(w_1=2, w_3=3,w_5=5\). 

    \begin{figure}[h]
        \centering
        \includegraphics[scale=1.4]{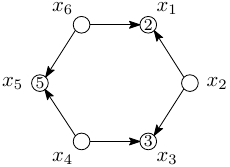}
    \end{figure}
    
    Then 
    $\mathbf{v}^{1}=(2,1,0,0,0,0)$, $\mathbf{v}^{2}=(0,1,3,0,0,0)$, $\mathbf{v}^{3}=(0,0,3,1,0,0)$, $\mathbf{v}^{4}=(0,0,0,1,5,0)$, $\mathbf{v}^{5}=(0,0,0,0,5,1)$, $\mathbf{v}^{6}=(2,0,0,0,0,1) \in \mathbb{N}^{6}$ are all lattice points of $P_{\calH}$. 
   By \Cref{claw} and \Cref{NonvertexK12}, there are no non-vertex lattice points of $P_{\calC_6}$.  By \Cref{CycelSinkSource},  $\calA[P_{\calC_6}]= K[P_{\calC_6}]$. 
\end{Example} 

\begin{acknowledgment*}
The authors thank Selvi Kara for helpful discussion at the beginning of this project. We thank the Institute for Computational and Experimental Research in Mathematics for their support and hospitality through the Collaborate@ICERM program in the summer of 2022.  
The authors used the software system Macaulay2 \cite{M2} to generate examples.
\end{acknowledgment*}

\noindent
\textbf{Data availability} Data sharing is not applicable to this article as no datasets were generated or analyzed
during the current study.

\noindent\textbf{Declaration of generative AI and AI-assisted technologies in the writing process:}
During the preparation of this work the authors used ChatGPT in order to improve readability and language of the work. After using this tool, the authors reviewed and edited the content as needed and take full responsibility for the content of the publication.

\bibliographystyle{plain}
\bibliography{references}

\end{document}